\documentclass[reqno,centertags, 12pt]{amsart}
\usepackage{amsmath,amsthm,amscd,amssymb}
\usepackage{latexsym}
\newcommand{\bbC}{{\mathbb{C}}}
\newcommand{\bbD}{{\mathbb{D}}}
\newcommand{\bbH}{{\mathbb{H}}}

\newcommand{\bbR}{{\mathbb{R}}}

\newcommand{\bbZ}{{\mathbb{Z}}}

\newcommand{\calI}{{\mathcal I}}

\newcommand{\lb}{\label}
\newcommand{\f}{\frac}

\newcommand{\ol}{\overline}
\newcommand{\ti}{\tilde  }

\newcommand{\tr}{\text{\rm{Tr}}}

\newcommand{\sgn}{\text{\rm{sgn}}}

\newcommand{\ess}{\text{\rm{ess}}}
\newcommand{\ac}{\text{\rm{ac}}}
\newcommand{\s}{\text{\rm{s}}}

\newcommand{\intt}{\text{\rm{int}}}
\newcommand{\ren}{\text{\rm{ren}}}
\newcommand{\supp}{\text{\rm{supp}}}

\newcommand{\bi}{\bibitem}

\newcommand{\beq}{\begin{equation}}
\newcommand{\eeq}{\end{equation}}
\newcommand{\ba}{\begin{align}}
\newcommand{\ea}{\end{align}}
\newcommand{\veps}{\varepsilon}

\newcounter{smalllist}
\newenvironment{SL}{\begin{list}{{\rm\roman{smalllist})}}{%
\setlength{\topsep}{0mm}\setlength{\parsep}{0mm}\setlength{\itemsep}{0mm}%
\setlength{\labelwidth}{2em}\setlength{\leftmargin}{2em}\usecounter{smalllist}%
}}{\end{list}}

\DeclareMathOperator{\Real}{Re}
\DeclareMathOperator{\Ima}{Im}

\let\det=\undefined\DeclareMathOperator{\det}{det}
\allowdisplaybreaks
\numberwithin{equation}{section}

\newtheorem{theorem}{Theorem}[section]
\newtheorem{proposition}[theorem]{Proposition}
\newtheorem{lemma}[theorem]{Lemma}
\newtheorem{corollary}[theorem]{Corollary}

\theoremstyle{remark}
\newtheorem*{remark}{Remark}
\newtheorem*{remarks}{Remarks}

\theoremstyle{definition}
\newtheorem*{definition}{Definition}

\newcommand{\abs}[1]{\lvert#1\rvert}

\begin{document}

\title[A Condition for Szeg\H{o} Asymptotics]
{Jost Functions and Jost Solutions\\for Jacobi Matrices,\\I.~A Necessary and Sufficient
Condition for Szeg\H{o} Asymptotics}
\author[D. Damanik and B. Simon]{David Damanik$^{1,2}$ and Barry Simon$^{1,3}$}

\thanks{$^1$ Mathematics 253-37, California Institute of Technology, Pasadena, CA 91125.
E-mail: damanik@caltech.edu; bsimon@caltech.edu}
\thanks{$^2$ Supported in part by NSF grant DMS-0227089}
\thanks{$^3$ Supported in part by NSF grant DMS-0140592 and in part by
Grant No.\ 2002068 from the United States-Israel Binational Science Foundation
(BSF), Jerusalem, Israel}

\date{March 28, 2003; February 10, 2005}

\begin{abstract} We provide necessary and sufficient conditions for a
Jacobi matrix to produce orthogonal polynomials with Szeg\H{o} asymptotics
off the real axis. A key idea is to prove the equivalence of Szeg\H{o}
asymptotics and of Jost asymptotics for the Jost solution. We also prove
$L^2$ convergence of Szeg\H{o} asymptotics on the spectrum.
\end{abstract}

\maketitle

\section{Introduction} \lb{s1}

In 1922, Szeg\H{o} \cite{Sz22a} proved one of the most celebrated results in classical
analysis: his asymptotic theorem for orthogonal polynomials. In modern language, he
considered measures, $d\rho$, on $[-2,2]$ of the form
\begin{equation} \lb{1.1}
d\rho(x) = f(x)\, dx + d\rho_\s (x)
\end{equation}
with orthonormal polynomials
\begin{equation} \lb{1.2}
p_n(x) = \gamma_n x^n + \text{ lower order}
\end{equation}
obeying $\gamma_n >0$ and
\begin{equation} \lb{1.3}
\int p_n(x) p_m(x)\, d\rho(x) = \delta_{nm}
\end{equation}
What Szeg\H{o} proved is that for $z\in\bbD =\{z\in\bbC\mid\abs{z} <1\}$, one has
Szeg\H{o} asymptotics as $n\to\infty$
\begin{equation} \lb{1.4}
z^n p_n \biggl( z + \f1{z}\biggr) \to \f{D(z)^{-1}}{\sqrt{2}}
\end{equation}
so long as the following, known as the Szeg\H{o} condition, holds
\begin{equation} \lb{1.5}
\int_{-2}^2 \log f(x) (4-x^2)^{-1/2} \, dx > -\infty
\end{equation}
(Actually, Szeg\H{o}, using the still standard convention of the orthogonal
polynomial community, took $d\rho$ on $[-1,1]$ and he did not allow a singular
component --- that is a later refinement. Also, instead of $z\mapsto z+z^{-1}$
which maps $\bbD\to \bbC\backslash [-2,2]$, he used the inverse map and stated his
results in terms of limits of
\begin{equation} \lb{1.6}
\biggl( \f{x}{2} + \f{\sqrt{4-x^2}}{2}\biggr)^n p_n(x)
\end{equation}
rather than \eqref{1.4}

Szeg\H{o} also found an explicit formula for $D(z)$, namely,
\begin{equation} \lb{1.7}
D(z) = \exp \biggl[ \int \f{e^{i\theta}+z}{e^{i\theta}-z}\, \log (f(\cos\theta))
\f{d\theta}{4\pi} \biggr]
\end{equation}
Moreover, if \eqref{1.5} fails, so does \eqref{1.4}

From the point of view of measures, the restriction to $\supp(d\rho)\subset
[-2,2]$ is natural, but this is less so with respect to the recursion coefficients
(aka Jacobi parameters) for the orthonormal polynomials, $p_n(x)$, defined by
\begin{equation} \lb{1.8}
xp_n(x) = a_{n+1} p_{n+1}(x) + b_{n+1} p_n(x) + a_n p_{n-1}(x)
\end{equation}
for $\{a_n,b_n\}_{n=1}^\infty$. From this point of view, the natural condition is
\begin{equation} \lb{1.9}
a_n\to 1 \qquad b_n\to 0
\end{equation}
This is associated to, indeed implies that, $\ess\,\supp(d\rho) = [-2,2]$, that is,
$\supp(d\rho) =[-2,2]\cup P$, where $P$ is a bounded set whose only possible
limit points are $\pm 2$. Our main goal in this paper is to answer the question
of for which $\{a_n,b_n\}_{n=1}^\infty$ does one have Szeg\H{o} asymptotics; we will
find (see Theorem~\ref{T5.1})

\begin{theorem}\lb{T1.1} Let $p_n(x)$ be orthonormal polynomials associated to Jacobi
parameters $\{a_n,b_n\}_{n=1}^\infty$ obeying \eqref{1.9}. Then $\lim z^n p_n (z+
\f{1}{z})$ exists for all $z\in\bbD$, is nonzero for $z\in\bbD\backslash\bbR$
with convergence uniform on compacts if and only if
\begin{alignat}{2}
&(\alpha) \qquad && \sum_{n=1}^\infty \, \abs{a_n -1}^2 + \abs{b_n}^2 <\infty \lb{1.10} \\
&(\beta) \qquad && \lim_{n\to\infty}\, a_n a_{n-1} \dots a_1 \text{ exists and is nonzero} \notag \\
&(\gamma) \qquad && \lim_{n\to\infty} \sum_{j=1}^n b_j \text{ exists} \notag
\end{alignat}
\end{theorem}
\noindent --- thereby closing a chapter opened 83 years ago.

There has, of course, been prior literature on these issues, although with considerably
stronger hypotheses than $(\alpha)$--$(\gamma)$. The initial results relating Jacobi
parameters to Szeg\H{o} asymptotics illustrated how strong $\supp(d\rho)\subset [-2,2]$
is and include

\begin{theorem} \lb{T1.2} Let $\supp (d\rho)\subset [-2,2]$. Then the following are
equivalent:
\begin{SL}
\item[{\rm{(a)}}] $(\beta)$ holds.
\item[{\rm{(b)}}] $(\alpha)$, $(\beta)$, and $(\gamma)$ hold.
\item[{\rm{(c)}}] The Szeg\H{o} condition \eqref{1.5} holds.
\end{SL}
\end{theorem}

This theorem combines results of Shohat \cite{Sh} and Nevai \cite{Nev79}; see also
\cite{KS} and \cite{SZ}. Of course, once one drops the restriction on $\supp(d\rho)$,
the $a$'s and $b$'s become almost independent, and any subset of $(\alpha)$--$(\gamma)$
can hold.

To continue our discussion of earlier results on extending Szeg\H{o} asymptotics,
we need some notation. Since $P$ can only have $\pm2$ as limit points,
\begin{equation} \lb{1.11}
P\cap (-\infty,-2) = \{E_j^-\}_{j=1}^{N_-}
\end{equation}
where $N_- =0$ (i.e., the set is empty), $1,2,\dots$ or $\infty$, and $E_1^- < E_2^-
<\cdots$. Similarly,
\begin{equation} \lb{1.12}
P\cap (2,\infty) = \{E_j^+\}_{j=1}^{N_+}
\end{equation}
with $E_1^+ > E_2^+ > \cdots$. The earliest results extending Szeg\H{o} asymptotics
beyond $\supp(d\rho)\subset [-2,2]$ are due to Gonchar \cite{Gon}, Nevai \cite{Nev79},
and Nikishin \cite{Nik}, who noted that the result still holds if $N_+ + N_- <\infty$.
More recently,

\begin{theorem}[Peherstorfer-Yuditskii \cite{PY}] \lb{T1.3} Suppose $a_n\to 1$,
$b_n\to 0$, and
\begin{equation} \lb{1.13}
\sum_{j,\pm}\, (\abs{E_j^\pm} - 2)^{1/2} <\infty
\end{equation}
and that \eqref{1.5} holds. Then \eqref{1.4} holds where the function $D(z)^{-1}$ vanishes
if and only if $z+z^{-1}$ is some $E_j^\pm$.
\end{theorem}

\begin{remark} The $D(z)^{-1}$ we use here is not the same as the $D^{-1}$ used in
\cite{PY}, but is a Blaschke product times their $D^{-1}$.
\end{remark}

Related to this is

\begin{theorem}[Killip-Simon \cite{KS}] \lb{T1.4} If
\begin{equation} \lb{1.14}
\sum_{n=1}^\infty\, \abs{a_n-1} + \abs{b_n} <\infty
\end{equation}
then \eqref{1.13} and \eqref{1.5} hold.
\end{theorem}

From one point of view, \eqref{1.13} is quite natural. If $z_j^\pm$ is defined by
\begin{equation} \lb{1.15}
z_j^\pm\in (-1,1) \qquad z_j^\pm + (z_j^\pm)^{-1} = E_j^\pm
\end{equation}
then \eqref{1.13} is equivalent to
\begin{equation} \lb{1.16}
\sum_{j,\pm} \, (1-\abs{z_j^\pm}) <\infty
\end{equation}
which is exactly what is needed to define a Blaschke product of zeros and obtain $D(z)^{-1}$
as a Nevanlinna function (see \cite{PY,KS,OPUC2}). Theorems~\ref{T1.3} and \ref{T1.4} are
the strongest prior results on when Szeg\H{o} asymptotics holds.

Both as input and motivation, the next element of background for our work concerns sum
rules. Szeg\H{o} proved his results for orthogonal polynomials on the real line (OPRL) by
mapping the problem to one on orthogonal polynomials on the unit circle (OPUC). For OPUC,
he earlier \cite{Sz20} proved asymptotic formulae. He began at $z=0$ where the limit
formula was equivalent to his leading limit theorem for Toeplitz determinants
(see \cite{Sz15}) and deduced the general formula from that.

Verblunsky \cite{V36} rewrote the $z=0$ limit theorem as a sum rule, namely, if
\begin{equation} \lb{1.17}
d\mu(\theta) = w(\theta)\, \f{d\theta}{2\pi} + d\mu_\s
\end{equation}
is a probability measure on $\partial\bbD$ and $\alpha_n$ are its Verblunsky coefficients
(see \cite{OPUC1,OPUC2} for definition), then
\begin{equation} \lb{1.18}
\prod_{j=0}^\infty \, (1-\abs{\alpha_j}^2) = \exp \biggl( \int \log(w(\theta))\,
\f{d\theta}{2\pi}\biggr)
\end{equation}
(which includes the fact that both sides are $0$ simultaneously, i.e., $\sum_{j=0}^\infty
\abs{\alpha_j}^2 = \infty \Leftrightarrow \int \log(w(\theta)) \f{d\theta}{2\pi} =
-\infty$). Without knowing of Verblunsky's work, Case \cite{Case1,Case2}, motivated by
KdV sum rules, wrote some rules for Jacobi matrices with sufficiently nice $a$'s and
$b$'s --- he was not explicit about the needed conditions, but his arguments at least
require
\begin{equation} \lb{1.19}
\sum_{n=1}^\infty n (\abs{a_n -1} + \abs{b_n}) <\infty
\end{equation}
It was Killip-Simon \cite{KS} who realized the right combination of sum rules and proved

\begin{theorem}[Killip-Simon \cite{KS}]\lb{T1.5} Let $a_n\to 1$ and $b_n \to 0$. Then
\begin{equation} \lb{1.20}
\sum_{n=1}^\infty \, \abs{a_n-1}^2 + \abs{b_n}^2 <\infty
\end{equation}
holds if and only if
\begin{equation} \lb{1.21}
\sum_{j,\pm} (\abs{E_j^\pm} -2)^{3/2} <\infty
\end{equation}
and
\begin{equation} \lb{1.22}
\int_{-2}^2 \log(f(x)) (4-x^2)^{1/2}\, dx <\infty
\end{equation}
\end{theorem}

Note that \eqref{1.22}, which \cite{KS} calls the quasi-Szeg\H{o} condition, is distinct
from \eqref{1.5} ($(4-x^2)^{1/2}$ rather than $(4-x^2)^{-1/2}$). Further developments of
sum rules include \cite{Kup1,Kup2,LNS,NPVY,Sim288,SZ}. In particular, one has

\begin{theorem}[Simon-Zlato\v{s} \cite{SZ}] \lb{T1.6} Consider the three assertions:
\begin{SL}
\item[{\rm{($\beta$)}}] $\lim_{n\to\infty} a_n \dots a_1$ exists and is nonzero.
\item[{\rm{($\sigma$)}}] \eqref{1.13} holds.
\item[{\rm{($\tau$)}}] \eqref{1.5} holds.
\end{SL}
If $(\beta)$ holds, then $(\sigma) \Leftrightarrow (\tau)$, and if $(\sigma)$
and $(\tau)$ hold, then $(\beta)$ holds.
\end{theorem}

The next element in our analysis is to link Szeg\H{o} asymptotics to a different
asymptotic result associated with work of Jost \cite{Jost}. Jost studied certain
solutions of $-u'' + Vu = Eu$, which is the analog of
\begin{equation} \lb{1.23}
a_n f_{n+1} + (b_n - (z+z^{-1})) f_n + a_{n-1} f_{n-1} =0
\qquad n=2,3,\dots
\end{equation}
one of whose solutions is
\begin{equation} \lb{1.24}
f_n (z) = p_{n-1} \biggl( z+ \f{1}{z}\biggr)
\end{equation}
As realized by Case \cite{Case1,Case2,GC80}, the analog of the Jost solution is a
solution of \eqref{1.23}, which is asymptotic to $z^n$ in the sense that
\begin{equation} \lb{1.25}
z^{-n} u_n(z) \to 1
\end{equation}
Case showed such solutions exist if $\abs{z} <1$ and \eqref{1.19} holds. In distinction,
Szeg\H{o} asymptotics says $p_{n-1} (z+\f{1}{z}) \sim Cz^{-n}$.

There may or may not be a solution of \eqref{1.23} which obeys \eqref{1.25} if one only
knows $a_n\to 1$, $b_n\to 0$, but from either the discrete version of Weyl's analysis
(see, e.g., \cite{RS2,S270}) or by the Poincar\'e-Perron theorem (see, e.g.,
\cite[Section~9.6]{OPUC2}), there is a solution for $z\in\bbD$ obeying $f_n\to 0$ ---
indeed, obeying $f_{n+1}/f_n\to z$. From Weyl's point of view, this is given by the
Green's function, that is, we can take it to be, for $z\in\bbD\backslash
\{z_j^\pm\}_{j=1}^{N_\pm}$,
\begin{equation} \lb{1.26}
w_n(z) = \langle \delta_n, (z+z^{-1} -J)\delta_1 \rangle
\end{equation}
where $J$ is the infinite Jacobi matrix
\begin{equation} \lb{1.27}
J=\begin{pmatrix} b_1 & a_1 & 0 & \dots \\
a_1 & b_2 & a_2 & \dots \\
0 & a_2 & b_3 & \dots \\
\dots & \dots & \dots & \dots \\
\dots & \dots & \dots & \dots
\end{pmatrix}
\end{equation}
viewed as a bounded selfadjoint operator on $\ell^2 (\bbZ_+)$.

We will say that Jost asymptotics occurs if for $z\in\bbD\backslash \{z_j^\pm\}_{j=1}^{N_\pm}$,
$z^{-n} w_n(z)$ has a nonzero finite limit as $n\to\infty$. A key to our understanding
of when Szeg\H{o} asymptotics holds for general $a$'s and $b$'s (i.e., to Theorem~\ref{T1.1})
is the following result we prove in Section~\ref{s2}:

\begin{theorem} \lb{T1.7} Fix $z_0\in\bbD$ so that $z_0 + z_0^{-1}$ is not an eigenvalue of $J$.
Then Szeg\H{o} asymptotics {\rm{(}}i.e., $z^n p_n (z+\f{1}{z})$ has a nonzero
limit{\rm{)}} holds at $z_0$ if and only if Jost asymptotics holds at $z_0$.
\end{theorem}

We can now turn more closely to our proof of Theorem~\ref{T1.1}. That Szeg\H{o} asymptotics
implies $(\alpha)$--$(\gamma)$ will be easy (and done in Section~\ref{s5}) once we have
Theorem~\ref{T1.7}. Basically, $\ti w_n (z)\equiv z^{-n} w_n(z)$ are analytic near
$z=0$ and Jost asymptotics (uniformly on $\abs{z}=\veps$) implies convergence of derivatives
at $z=0$. The first two Taylor terms at $0$ yield $(\beta)$--$(\gamma)$, and as in \cite{KS}, a
suitable combination of the first and third Taylor coefficients is positive and yields
$(\alpha)$.

The hard direction is that $(\alpha)$--$(\gamma)$ implies Szeg\H{o} or Jost asymptotics.
We will provide three distinct proofs. The first, in Section~\ref{s5}, is a relative of
Szeg\H{o}'s original proof and of the Peherstorfer-Yuditskii arguments relying on the study
of analytic functions on the disk. Szeg\H{o} just used \eqref{1.7} to define $D$, and
Peherstorfer-Yuditskii multiplied $D^{-1}$ by a Blaschke product. We do not have either
luxury here. For \eqref{1.7} to work, one needs
\begin{equation} \lb{1.28}
\int \log(f(\cos\theta))\, \f{d\theta}{2\pi} >-\infty
\end{equation}
which is equivalent to \eqref{1.5}, while all we have is
\begin{equation} \lb{1.29}
\int \log(f(\cos\theta)) \sin^2 (\theta)\, \f{d\theta}{2\pi} >-\infty
\end{equation}
which is equivalent to \eqref{1.22}. Moreover, in place of \eqref{1.16}, we only have
\begin{equation} \lb{1.30}
\sum_{j,\pm}\, (1-\abs{z_j})^3 <\infty
\end{equation}
so we cannot define a Blaschke product. The solution will be to define renormalized
Blaschke products when \eqref{1.30} holds, which we do in Section~\ref{s3}, and a
renormalized Poisson integral when \eqref{1.29} holds, which we do in Section~\ref{s4}.
This will allow us to define a candidate for the Jost function and prove Jost asymptotics
in Section~\ref{s5} and so provide our first proof that $(\alpha)$--$(\gamma)$ imply
Jost asymptotics. This proof provides bounds we will need in Section~\ref{s8} to handle
$L^2$ convergence on $\partial\bbD$.

Our second proof in Section~\ref{s6} relies on an idea going back to Jost-Pais
\cite{JP} that the Jost function is a Fredholm determinant. For OPRL, this is discussed
in Killip-Simon \cite{KS}. We will use the theory of renormalized determinants for
Hilbert-Schmidt operators to construct a candidate Jost function and use it to prove
Jost asymptotics.

Our final proof, in Section~\ref{s7}, is connected to classical results on the construction
of asymptotic solutions of ODE's associated with work of Levinson \cite{Lev} and
Hartman-Wintner \cite{HW}; see the book of Eastham \cite{East}. We will use results of
Coffman \cite{Coff} on the difference equation analogs to construct Jost solutions when
$(\alpha)$--$(\gamma)$ hold. This construction shows that the ``hard" part of
Theorem~\ref{T1.1} is related to known results on ODE's with $L^2$ perturbations. From
this point of view, our contribution here is the realization that Jost solutions
imply Szeg\H{o} asymptotics and that the conditions are not only sufficient but necessary.

In Section~\ref{s8}, we discuss $L^2$ convergence on $\partial\bbD$, following the
original scheme of Szeg\H{o} \cite{Sz20} but with some severe technical complications
because the Jost function is not Nevanlinna. This is the hardest argument in the paper.
In Section~\ref{s9}, we provide examples for each $p<\f32$ of Jacobi matrices with
Szeg\H{o} asymptotics, but with $\sum_{j,\pm} (\abs{E_j^\pm}-2)^p =\infty$.
In Section~\ref{s10}, we make some remarks about Schr\"odinger operators with $L^2$
potentials.

We announced our results in \cite{DKS} written in September of 2003 and mentioned our
$L^2$ results but not their proof to Serguei Denisov. In May of 2004, Denisov-Kupin
\cite{DeK} released a preprint discussing modified Szeg\H{o} asymptotics for certain OPUC
when the Szeg\H{o} condition fails but a condition like \eqref{1.29} holds. Their results
are quite distinct from ours although, via \eqref{1.29}, there is some overlap. Many of
the methods are similar
--- in particular, like we do in Section~\ref{s4}, they use renormalized Poisson
representations. There is also some overlap in the $L^2$ control of the boundary values
which we consider in Section~\ref{s8}. In particular, by using some of their ideas, it is
likely we could streamline the proof of and slightly strengthen our estimate,
Proposition~\ref{P8.2}. We have kept our original proof. We would emphasize that our work
on these methods is independent and roughly simultaneous.

\smallskip
It is a pleasure to thank M.~Moszy\'nski and R.~Romanov for useful discussions. B.S.\
completed this work during his stay as a Lady Davis Visiting Professor at Hebrew
University, Jerusalem. He would like to thank H.~Farkas and Y.~Last for the hospitality
of the Mathematics Institute at Hebrew University.

\section{Szeg\H{o} Asymptotics and Jost Asymptotics} \lb{s2}

As explained in the introduction, for any Jacobi matrix with $a_n\to 1$, $b_n\to 0$, and $z\in\bbD$,
and not such that $z+z^{-1}$ is an eigenvalue of $J$, there are two natural solutions of
\begin{equation} \lb{2.1}
a_n f_{n+1} + (b_n -(z+z^{-1}))f_n +a_{n-1} f_{n-1} =0 \qquad n=2,3,\dots
\end{equation}
One is the orthogonal polynomial solution, $f_n=p_{n-1} (z+\f{1}{z})$, and the other is the Weyl solution,
\begin{equation} \lb{2.2}
w_n(z) =\langle\delta_n, (z+z^{-1} -J)^{-1} \delta_1 \rangle
\end{equation}

In this section, our purpose is to show that for each such $z$, one has Jost asymptotics at that $z$,
that is
\begin{equation} \lb{2.3}
\ti w_n(z) \equiv z^{-n} w_n (z) \to \ti w_\infty (z)
\end{equation}
for $\ti w_\infty \neq 0$ if and only if one has Szeg\H{o} asymptotics for that $z$, that is,
\begin{equation} \lb{2.4}
c_n (z) \equiv z^n p_n \biggl(z+\f{1}{z}\biggr)\to c_\infty (z)
\end{equation}
for $c_\infty \neq 0$, and moreover,
\begin{equation} \lb{2.5}
(1-z^2) c_\infty (z) \ti w_\infty (z)=1
\end{equation}
(as we will see, $\ti w_\infty (z) =1/u(z)$, where $u$ is the Jost function, so \eqref{2.5} is usually
written $c_\infty (z) =u(z)/(1-z^2)$).

Of course, $p_{\boldsymbol{\cdot}\, -1}$ obeys \eqref{2.1} also at $n=1$ if we define $p_{-1}\equiv f_0
=0$ and $a_0=1$. Since
\[
(J-z-z^{-1})(z+z^{-1} -J) \delta_1 =-\delta_1
\]
$w_n$ also obeys \eqref{2.1} if we set $a_0=1$ and
\begin{equation} \lb{2.5a}
w_0(z) =1
\end{equation}

The constancy of the Wronskian thus implies
\begin{equation} \lb{2.6}
a_n \biggl(p_n \biggl(z+\f{1}{z}\biggr) w_n(z) -w_{n+1}(z) p_{n-1} \biggl( z+\f{1}{z}\biggr)\biggr) =1
\end{equation}
where we get $1$ since
\[
a_0 (p_0 w_0 - w_1 p_{-1})=1
\]
Using the definitions \eqref{2.3}/\eqref{2.4} of $c$ and $\ti w$, \eqref{2.6} becomes
\begin{equation} \lb{2.7}
a_n (c_n(z) \ti w_n (z) -z^2 \ti w_{n+1}(z) c_{n-1} (z))=1
\end{equation}

Thus, the following lemma is of relevance:

\begin{lemma} \lb{L2.1} Let $x_n,y_n$ be sequences of nonzero complex numbers and let $\lambda_n$ be
nonzero positive numbers with
\begin{equation} \lb{2.8}
\lambda_n\to 1
\end{equation}
and so, for some $z\in\bbD$,
\begin{equation} \lb{2.9}
x_{n+1} y_n -z^2 x_n y_{n+1} =\lambda_n
\end{equation}
Then
\begin{SL}
\item[{\rm{(i)}}] If $y_n\to y_\infty \neq 0$, then $x_n \to 1/y_\infty (1-z^2)$.
\item[{\rm{(ii)}}] If $x_n\to x_\infty \neq 0$ and $z^{2n} y_n\to 0$, then $y_n \to 1/x_\infty (1-z^2)$.
\end{SL}
\end{lemma}

\begin{proof} (i) Rewrite \eqref{2.9} as
\[
x_{n+1} = \lambda_n y_n^{-1} + z^2 \, \f{y_{n+1}}{y_n}\, x_n
\]
and iterate $\ell+1$ times to get
\begin{equation} \lb{2.10}
x_{n+1} =\sum_{j=0}^\ell \lambda_{n-j}\, \f{y_{n+1}}{y_{n+1-j} y_{n-j}}\, z^{2j}
+ z^{2\ell +2} \, \f{y_{n+1}}{y_{n-\ell}}\, x_{n-\ell}
\end{equation}
Set $\ell=n-1$ and see that since $\sum_{j=0}^{n-1} z^{2j} = (1-z^{2n})/(1-z^2)$,
\begin{equation} \lb{2.11}
\abs{x_{n+1} -y_\infty^{-1} (1-z^{2n})(1-z^2)^{-1}} \leq \sum_{j=0}^{n-1} e_{n,j} z^{2j}
+ e_{n,n} z^{2n}
\end{equation}
where
\begin{align*}
e_{n,j} & =\lambda_{n-j} \, \f{y_{n+1}}{y_{n+1-j} y_{n-j}} - \f{1}{y_\infty} \qquad j=0, \dots, n-1 \\
e_{n,n} &= y_{n+1} x_1 y_1^{-1}
\end{align*}
Since $y_\ell\to y_\infty \neq 0$, $\sup_{n,j} e_{n,j}<\infty$ and moreover,
$\lim_{n\to\infty} e_{n,j} =0$ for all fixed $j$. Thus, since $e_{n,j}\to 0$
for $j$ fixed, we have for $\ell$ fixed,
\begin{align*}
\limsup_{n\to\infty} \biggl|\, \sum_{j=0}^n e_{n,j} z^{2j} \biggr|
&\leq \limsup \, \biggl| \, \sum_{j=\ell}^n e_{n,j} z^{2j}\biggr| \\
&\leq \abs{z}^{2\ell} (1-\abs{z}^2)^{-1} \sup_{n,j}\, \abs{e_{n,j}}
\end{align*}
$\to 0$ as $\ell\to\infty$. Thus \eqref{2.11} implies
\begin{equation} \lb{2.12}
x_{n+1}\to y_\infty^{-1} (1-z^2)^{-1}
\end{equation}

\smallskip
(ii) Rewrite \eqref{2.9} as
\[
y_n =\lambda_n x_{n+1}^{-1} + z^2 x_n x_{n-1}^{-1} y_{n+1}
\]
and iterate upwards. Since $z^{2n} y_n\to 0$, the remainder after $\ell$
iterations goes to zero as $\ell\to\infty$, so
\[
y_n =\sum_{j=0}^\infty \lambda_{n+j} z^{2j} x_n x_{n+j+1}^{-1} x_{n+j}^{-1}
\]
As in the argument in (i), this implies that $y_n\to x_\infty^{-1} (1-z^2)^{-1}$.
\end{proof}

\begin{theorem}[Szeg\H{o} asymptotics $=$ Jost asymptotics]\lb{T2.2} Let $J$ be
a Jacobi matrix with $a_n\to 1$, $b_n\to 0$, and let $z\in\bbD$ be such that
$z+z^{-1}$ is not an eigenvalue of $J$. Then $\ti w_n (z)$ has a nonzero limit
if and only if $c_n(z)$ has a nonzero limit, and if either happens,
\begin{equation} \lb{2.13}
\lim_{n\to\infty} \, c_n(z) = \f{u(z)}{1-z^2}
\end{equation}
where
\begin{equation} \lb{2.14}
u(z)^{-1} \equiv\lim_{n\to\infty} \, \ti w_n(z)
\end{equation}
\end{theorem}

\begin{proof} By \eqref{2.7}, if $\lambda_n =a_n^{-1}$, $x_n =c_{n-1} (z)$,
$y_n =\ti w_n(z)$, then Lemma~\ref{L2.1} implies this result so long as
\[
\lim_{n\to\infty} \, z^{2n} \ti w_n (z) =0
\]
But
\[
z^{2n} \ti w_n (z) =z^n w_n(z)
\]
goes to zero since both $w_n\to 0$ and $z^n \to 0$.
\end{proof}

\section{Renormalized Blaschke Products} \lb{s3}

As explained in the introduction, we need a renormalized Blaschke product
that works for real zeros that only obey $\sum_n (1-\abs{z_n})^3 <\infty$
rather than the usual Blaschke condition $\sum_n (1-\abs{z_n})<\infty$.

One can make a case that the first renormalization in science was the
Weierstrass product formula --- to get an analytic function vanishing
at $\{z_j\}_{j=1}^\infty$ with $\abs{z_j}\to\infty$, one modifies
one's first guess
\[
f(z)=\prod_{j=1}^\infty \biggl( 1-\f{z}{z_j}\biggr)
\]
to
\begin{equation} \lb{3.1}
f(z) =\prod_{j=1}^\infty W_{n_j} \biggl( \f{z}{z_j}\biggr)
\end{equation}
where
\begin{equation} \lb{3.2}
W_n(z) =(1-z) \exp \biggl(\, \sum_{k=1}^n \f{z^k}{k}\biggr)
\end{equation}
picking the argument to be the truncation of the power series for $-\log (1-z)$. It is
well known, of course, that \eqref{3.1} converges if $n_j$ is chosen so that $\sum
\abs{r/z_j}^{n_j+1}<\infty$ for all $r>0$. Similarly, if our only goal were to get a
function with zeros in the right place, things would be easy --- for one can show that if
$z_j\in\bbD$, $\abs{z_j}\to 1$ as $j\to \infty$ and $w_j=z_j/\abs{z_j}$, and if $n_j$ is
chosen so that $\sum_{j=1}^\infty (\f{1-\abs{z_j}}{\veps})^{n_j+1}<\infty$ for all $\veps
>0$, then
\begin{equation} \lb{3.3}
f(z) =\prod_{j=1}^\infty W_{n_j} \biggl( \f{w_j -z_j}{w_j-z}\biggr)
\end{equation}
is a product converging absolutely to a nonzero function analytic in
$\bbD$ with zeros at $\{z_j\}$.

We want our Blaschke products to have magnitude one on $\partial\bbD$ and we
will want that for our renormalized Blaschke products.

For $p=0$, $b(z,p)=z$: if $p\in\bbD$, $p\neq 0$,
\begin{equation} \lb{3.4}
b(z,p) =\f{\abs{p}}{p}\, \f{p-z}{1-\bar pz}
\end{equation}
so $b(0,p)=\abs{p}>0$. The key, of course, is that $b(z,p)=0$ if and only if
$z=p$ and
\begin{equation} \lb{3.5}
\abs{b(e^{i\theta},p)}=1
\end{equation}

If $p=(1-x)\omega$ with $\abs{\omega} =1$ and $x\in (0,1)$,
\begin{align}
b(z, (1-x)\omega) &= \f{1-x-z\omega^{-1}}{1-(1-x)\omega^{-1} z} \notag \\
&= \f{1-\f{x}{1-\omega^{-1}z}}{1+\f{x\omega^{-1}z}{1-\omega^{-1}z}} \lb{3.6}
\end{align}
\eqref{3.6} shows immediately if $\abs{z}<1$ and $\sum \abs{x_j} <\infty$,
then $\prod_{j=1}^\infty b(z, (1-x_j)\omega_j)$ converges absolutely
(and uniformly on $\abs{z}<1 -\delta$) since the numerators and
denominators in \eqref{3.5} separately do.

\eqref{3.6} suggests what to do if $\sum \abs{x_j}^{n+1}<\infty$. Define
\begin{equation} \lb{3.8}
b_n (z,(1-x)\omega) =
\f{W_n (\f{x}{1-\omega^{-1}z})}{W_n (\f{-x\omega^{-1}z}{1-\omega^{-1} z})}
\end{equation}

Here is the key fact:

\begin{proposition} \lb{P3.1}
\begin{SL}
\item[{\rm{(a)}}] Let $\delta >0$. Then for $\abs{z}<1-\delta$ and $\abs{x}
<\delta/2$,
\begin{equation} \lb{3.9}
\abs{b_n (z, (1-x)\omega) -1} \leq 4\delta^{-n-1} x^{n+1}
\end{equation}
\item[{\rm{(b)}}] For $e^{i\theta}\neq\omega$,
\begin{equation} \lb{3.10}
\abs{b_n (e^{i\theta}, (1-x)\omega)}=1
\end{equation}
\end{SL}
\end{proposition}

{\it Warning.} One cannot use the maximum principle and \eqref{3.10} to conclude
that $\abs{b_n (z, (1-x)\omega)}\leq 1$. Indeed, for $n\geq 1$,
\[
\lim_{r\uparrow 1}\, \abs{b_n (r\omega, (1-x)\omega)}=\infty
\]
This is where $b_n$'s differ from ordinary Blaschke factors. They have very
singular inner factors (indeed, for $n\geq 3$, ones whose boundary values are
not even signed measures).

\begin{proof} (a) It is known, (e.g., Rudin \cite[p.~301]{Rudin}) that
\begin{equation} \lb{3.11}
\abs{z}< 1 \Rightarrow \abs{W_n(z) -1} \leq \abs{z}^{n+1}
\end{equation}
If $\abs{x}<\delta/2$ and $\abs{z} <1-\delta$, then $\abs{x/(1-\omega^{-1} z)}\leq \abs{x/\delta}
<\f12$ so \eqref{3.11} can be used, and if $N$ and $D$ are the numerator and denominator in \eqref{3.8},
$\abs{D}>\f12$. Since
\[
\biggl| \f{N}{D}-1\biggr| \leq \f{1}{\abs{D}} \, (\abs{N-1} + \abs{D-1}) \]
\eqref{3.11} implies that
\[
\abs{b_n (z, (1-x)\omega)} \leq 2 \biggl[ \, \biggl| \f{x}{(1-\omega^{-1} z)}\biggr|^{n+1} +
\biggl| \f{x\omega^{-1} z}{(1-\omega^{-1} z)}\biggr|^{n+1} \biggr]
\]
which yields \eqref{3.9}.

\smallskip
(b) By \eqref{3.6}, if $e^{i\theta}\neq \omega$, $b(e^{i\theta}, (1-x)\omega)$ can be defined as a limit
and \eqref{3.6} still holds and $b(e^{i\theta},\omega)=1$. Thus for $x$ small, $\log b(e^{i\theta},
(1-x)\omega)$ is analytic in $x$. By \eqref{3.4},
\[
\abs{b(e^{i\theta}, p)} = \biggl| \f{p-e^{i\theta}}{1-\bar pe^{i\theta}}\biggr| =
\biggl| \f{p-e^{i\theta}}{\bar p-e^{-i\theta}}\biggr| =1
\]
so for $x$ positive, $e^{i\theta}\neq \omega$,
\[
\Real \log (b(e^{i\theta}, (1-x)\omega))=0
\]
It follows that its Taylor coefficients,
\begin{equation} \lb{3.12}
\log (b(e^{i\theta}, (1-x)\omega)) =\sum_{n=1}^\infty \gamma_n (e^{i\theta}, \omega) x^n
\end{equation}
have $\gamma_n$ pure imaginary. Since
\begin{equation} \lb{3.13}
b_n (e^{i\theta}, (1-x)\omega)\equiv b(e^{i\theta}, (1-x)\omega) \exp \biggl( \, \sum_{j=1}^n
\gamma_j (e^{i\theta}, \theta) x^j \biggr)
\end{equation}
and $\gamma_j$ is pure imaginary, \eqref{3.10} holds.
\end{proof}

Because we will be interested not in $b_2$ but something related to it by a finite correction, we need
to look in detail at $\gamma_1$ and $\gamma_2$. We consider $\gamma_j (z,\omega)$ defined by \eqref{3.12}
with $e^{i\theta}\to z$. By \eqref{3.6},
\begin{align}
\gamma_1 (z,\omega) &= -\biggl( \f{1+\omega^{-1}z}{1-\omega^{-1} z}\biggr)  \lb{3.14} \\
\gamma_2 (z,\omega) &= -\tfrac12\, \f{(1-(\omega^{-1} z)^2)}{(1-\omega^{-1} z)^2} = -\tfrac12 \,
\f{(1+\omega^{-1} z)}{(1-\omega^{-1}z)} \lb{3.15}
\end{align}
Remarkably, $\gamma_1/\gamma_2$ is independent of $\omega$ and $z$! For reasons that will be clear
below, we want to consider
\begin{equation} \lb{3.16}
\alpha(z) = \f{1+z^2}{1-z^2} \qquad \beta(z) =\f{2z}{1-z^2}
\end{equation}
Notice that
\begin{equation} \lb{3.17}
\gamma_1 (z,\omega =\pm 1) = 2\gamma_2 (z,\omega =\pm 1) = -(\alpha (z) \pm \beta(z))
\end{equation}

\begin{definition} For $p\in (-1,1)$, $p\neq 0$, and $z\in\bbD$, we define
\begin{equation} \lb{3.18}
q(z,p) =b(z,p) \exp (-\alpha (z) \log (\abs{p})-\tfrac12\, \beta(z)(p-\tfrac{1}{p}))
\end{equation}
\end{definition}

\begin{theorem} \lb{T3.2}
\begin{SL}
\item[{\rm{(a)}}] For $z$ near zero and $p\neq 0$, $p$ real,
\begin{equation} \lb{3.19}
\log q(z,p) =\log b(z,p) - \alpha(z) \log b(0,p) -\tfrac12\, \beta(z) \left. \f{d}{dz}\,
\log b(z,p)\right|_{z=0}
\end{equation}
\item[{\rm{(b)}}]
\begin{equation} \lb{3.20}
q(z,p) =b_2 (z,p) \exp (-\alpha(z) A(p) -\tfrac12\, \beta(z) B(p))
\end{equation}
where
\begin{align}
A(p) &= \log \abs{p} - (1-\abs{p}) - \f{(1-\abs{p})^2}{2} \lb{3.21}  \\
B(p) &= -\f{(1-\abs{p})^3}{p} \lb{3.22}
\end{align}
\item[{\rm{(c)}}] If $p\in (-1,1)$ with $1-\abs{p}<\delta/2$ and $\abs{z}< 1-\delta$, then
\begin{equation} \lb{3.23}
\abs{q(z,p)-1} \leq [4\delta^{-3} + \tfrac53(1+4\delta^{-3}) \delta^{-1} \abs{p}^{-1}
\exp (\tfrac53\, \delta^{-1} \abs{p}^{-1} (1-\abs{p})^3] (1-\abs{p})^3
\end{equation}
\end{SL}
\end{theorem}

\begin{proof} (a) Writing
\[
b(z,p) = \abs{p} \, \f{1-\f{z}{p}}{1-zp}
\]
we see
\[
\log b(z,p) = \log \abs{p} + z\biggl( p-\f{1}{p}\biggr) + O(z^2)
\]
which, given \eqref{3.18}, is \eqref{3.19}.

\smallskip
(b) By \eqref{3.13} and \eqref{3.17},
\[
q(z,p) =b_2 (z,x) \exp (C(p,z))
\]
where
\begin{equation} \lb{3.24}
\begin{split}
C(p,z) & =-\alpha(z) \log (\abs{p}) -\tfrac12\, \beta(z) \biggl( p-\f{1}{p}\biggr)
 -\alpha(z) ((1-\abs{p}) + \tfrac12 (1-\abs{p})^2) \\
 &\qquad - \beta(z) \sgn(p) ((1-\abs{p}) + \tfrac12 (1-\abs{p})^2)
\end{split}
\end{equation}
Thus \eqref{3.20} follows from
\begin{equation} \lb{3.25}
p-\f{1}{p} + \sgn\abs{p} \, [2(1-\abs{p}) + (1-\abs{p})^2] = -\f{(1-\abs{p})^3}{p}
\end{equation}
\eqref{3.24} follows from writing $p=\sgn (p) (1-x)$ and
\begin{align*}
\sgn(p) \biggl[ (1-x) &- \f{1}{1-x} + 2x + x^2\biggr] \\
&= \sgn(p) \biggl\{ (1-x) - \biggl[ 1+x+x^2 + \f{x^3}{1-x}\biggr] + 2x + x^2\biggr\} \\
&= -\sgn(p) \, \f{x^3}{1-x} = -\f{(1-\abs{p})^3}{p}
\end{align*}

\smallskip
(c) In terms of the function $C$ of \eqref{3.24},
\begin{align}
\abs{q(z,p)-1} &= \abs{b_2 (z,p) \exp (C(z,p)) -1} \\
&\leq \abs{b_2 (z,p)} \, \abs{\exp (C(p,z))-1} + \abs{b_2 (z,p)-1} \lb{3.25a}
\end{align}
We have \eqref{3.9} to bound $\abs{b_2 (z,p)-1}$. Thus $\abs{b_2} \leq 1+\abs{b_2-1}
\leq 1 + 4\delta^{-3}$. Moreover,
\[
\abs{e^c-1} \leq \abs{c} \max (1,\abs{e^c}) \leq \abs{c} e^{\abs{c}}
\]
so \eqref{3.23} follows from
\begin{equation} \lb{3.26}
\abs{C} \leq \tfrac53\, \delta^{-1} \abs{p}^{-1} (1-\abs{p})^3
\end{equation}

To prove \eqref{3.26}, note first that $\abs{1-z^2} \geq 1-\abs{z}^2=(1+\abs{z})(1-\abs{z})\geq\delta$.
Thus
\begin{equation} \lb{3.27}
\abs{\alpha(z)} \leq \f{2}{\delta} \qquad \abs{\beta(z)} \leq \f{2}{\delta}
\end{equation}
Moreover, if $\abs{p}=1-x$, then
\begin{align*}
\biggl|\log \abs{p} -x - \f{x^2}{2}\biggr| &= \biggl|\, \sum_{j=3}^\infty \f{x^j}{j}\biggr| \\
&\leq \f13\, \f{x^3}{1-x} \\
&= \f13\, \f{(1-\abs{p})^3}{\abs{p}}
\end{align*}
Thus, by \eqref{3.25} and \eqref{3.27},
\[
\abs{C} \leq \f{2}{\delta}\, \f{(1-\abs{p})^3}{p} \, \biggl[\f13 + \f12 \biggr] = \f{5}{3\delta} \,
\f{(1-\abs{p})^3}{\abs{p}}
\]
proving \eqref{3.26}.
\end{proof}

Because each $b_n (z,p)$ is unbounded on $\bbD$, the usual methods for
controlling products on $\partial\bbD$ do not work; but in the case where
the limit points of zeros only are a finite set, they do. Here is what
we will need:

\begin{theorem} \lb{T3.3} Let $p_n$ be a sequence of reals in $(-1,1)$
with $\lim_{n\to\infty} \abs{p_n}=1$ so that
\begin{equation} \lb{3.28}
\sum_{n=1}^\infty (1-\abs{p_n})^3 <\infty
\end{equation}
Let
\begin{equation} \lb{3.29}
B_{\ren}(z) =\prod_{n=1}^\infty q(z,p_n)
\end{equation}
Then
\begin{SL}
\item[{\rm{(i)}}] The product \eqref{3.29} converges in $\bbC_+ =\{z\in \Ima z>0\}$
and defines a function analytic in $\bbD \cup \bbC_+\cup\bbC_-$ whose only zeros
are at $\{p_n\}_{n=1}^\infty$.
\item[{\rm{(ii)}}]
\begin{equation} \lb{3.30}
\abs{B_\ren (e^{i\theta})}=1 \qquad \theta\in (0,\pi)\cup (\pi, 2\pi)
\end{equation}
\item[{\rm{(iii)}}] If
\begin{equation} \lb{3.31}
B_{\ren}^{(N)} =\prod_{n=N+1}^\infty q(z,p_n)
\end{equation}
then for any $p<\infty$,
\begin{equation} \lb{3.32}
\int \abs{B_{\ren}^{(N)} (e^{i\theta}) -1}^p \, \f{d\theta}{2\pi} \to 0
\end{equation}
\end{SL}
\end{theorem}

\begin{proof} (i) If $z\in\bbC_+$, we have
\[
G(z) \equiv \max \biggl( \f{1}{\abs{1-z}}\, , \f{1}{\abs{1+z}}\, ,
\f{\abs{z}}{\abs{1-z}}\, , \f{\abs{z}}{\abs{1+z}} \biggr) <\infty
\]
Thus, if $xG(z)<1$, the arguments in $W_n$ in \eqref{3.8} are less than $1$ and
the same estimates we used to bound $\abs{q(z,p)-1}$ still work to see
\begin{equation} \lb{3.33}
\abs{q(z,p) -1} \equiv H(z) \abs{1-p}^3
\end{equation}
for suitable $H(z)$, and this shows the product converges.

\smallskip
(ii) Since the product converges on $\partial\bbD\backslash \{\pm 1\}$ and
$\abs{q(e^{i\theta},p_n)} =1$, \eqref{3.30} is immediate.

\smallskip
(iii) Since $\abs{B_{\ren}^{(N)} (e^{i\theta})}=1$, by (ii), pointwise convergence
implies $L^p$ convergence. The estimate \eqref{3.33} implies pointwise convergence
to $1$ since $\sum_{n=N+1}^\infty \abs{q(z,p_n) -1} \to 0$.
\end{proof}

\section{Renormalized Poisson Representations} \lb{s4}

Our goal in this section is to start out with a function, $f(z)$, on $\bbD$,
which has a complex Poisson representation
\begin{equation} \lb{4.1}
f(z) =\int P(z,e^{i\theta}) g(e^{i\theta}) \, \f{d\theta}{2\pi}
\end{equation}
where
\begin{equation} \lb{4.2}
P(z, e^{i\theta}) = \f{e^{i\theta} + z}{e^{i\theta} -z}
\end{equation}
and $g\in L^1 (d\theta/2\pi)$, real-valued, and
\begin{equation} \lb{4.3}
g(e^{i\theta}) = g(e^{-i\theta})
\end{equation}
(so $f(z)$ is real on $\bbD\cap\bbR$).

We want to define
\begin{equation} \lb{4.4}
h(z)=f(z)-\alpha(z) f(0) -\beta (z) f'(0)
\end{equation}
and show it has a representation
\begin{equation} \lb{4.5}
h(z) =\int Q(z,e^{i\theta}) g(e^{i\theta}) \, \f{d\theta}{2\pi}
\end{equation}
where $Q$ obeys a bound
\begin{equation} \lb{4.6}
\abs{Q(z,e^{i\theta})}\leq C(z) \sin^2 \theta
\end{equation}
This will allow us to extend \eqref{4.1} to cases where one only has $\int
\abs{g(e^{i\theta})}\sin^2 \theta \f{d\theta}{2\pi} <\infty$. In \eqref{4.4}, $\alpha$
and $\beta$ are the functions in \eqref{3.16}. For this section, their key property is
\begin{equation} \lb{4.7}
\alpha(z) \pm \beta(z) = P(z,e^{i\theta} =\pm 1)
\end{equation}

To see why \eqref{4.6} should hold, note that, by \eqref{4.3}, in \eqref{4.1}
we can replace $P(z,e^{i\theta})$ by
\begin{align}
S(z, e^{i\theta}) &\equiv \tfrac12\, [P(z, e^{i\theta}) + P(z, e^{-i\theta})]  \lb{4.8} \\
&= \f{1-z^2}{1+z^2 -2z \cos\theta} \lb{4.9}
\end{align}
Since $S(0, e^{i\theta}) =1$ and $\left.\f{\partial}{\partial z}
S(z,e^{i\theta})\right|_{z=0} = \cos\theta$, \eqref{4.5} holds with
\begin{equation} \lb{4.10}
Q(z, e^{i\theta}) =S(z, e^{i\theta}) -\alpha(z) -\beta(z) \cos\theta
\end{equation}
Because of \eqref{4.7} and $P(z, e^{i\theta}=\pm 1)=S(z, e^{i\theta} =
\pm 1)$, $Q$ vanishes at $e^{i\theta} =+ 1$ and at $e^{i\theta} =-1$.
Since $\alpha$ is even under $\theta\to -\theta$ and $\theta\to 2\pi -\theta$,
these zeros must be quadratic, which is where \eqref{4.6} comes from.

A straightforward calculation shows that, by \eqref{4.10},
\begin{align}
Q(z, e^{i\theta}) &= S(z, e^{i\theta}) - \f{1+z^2}{1-z^2} - \f{2z\cos\theta}{1-z^2} \notag \\
&= \f{1-z^2}{1+z^2 - 2z\cos\theta} - \f{1+z^2 + 2z \cos\theta}{1-z^2} \notag \\
&= \f{-4z^2 \sin^2 \theta}{(1-z^2) (1+z^2 -2z\cos\theta)}  \lb{4.11} \\
&= \f{-4z^2 \sin^2 \theta}{(1-z)(1+z)(z-e^{i\theta}) (z-e^{-i\theta})} \lb{4.12}
\end{align}
We summarize with

\begin{theorem}\lb{T4.1} Let $f$ be given by \eqref{4.1} with $g\in L^1
(d\theta/2\pi)$ satisfying \eqref{4.3} and let $h(z)$ be given by \eqref{4.4}. Then
\eqref{4.5} holds with $Q$ given by \eqref{4.11}. In particular,
\begin{equation} \lb{4.13}
\abs{Q(z, e^{i\theta})} \leq \f{4\sin^2\theta}{(1-\abs{z})^3}
\end{equation}
\end{theorem}

\begin{proof} To get \eqref{4.13}, note that $\abs{1-z^2} \geq 1-\abs{z}^2
=(1-\abs{z})(1+\abs{z}) \geq 1-\abs{z}$ and $\abs{z-e^{\pm i\theta}} \geq 1-\abs{z}$.
\end{proof}

As a final result about renormalized Poisson representations, we note that

\begin{theorem} \lb{T4.2} Let $g\in L^1 (\sin^2\theta\, \f{d\theta}{2\pi})$ be
real-valued with $g(e^{i\theta}) =g(e^{-i\theta})$. Define
\begin{equation} \lb{4.14}
f(z) =\int_0^{2\pi} Q(z, e^{i\theta}) g(e^{i\theta}) \, \f{d\theta}{2\pi}
\end{equation}
Then for a.e.~$\theta$, $\lim_{r\uparrow 1} f(re^{i\theta})\equiv f(e^{i\theta})$
exists, and for a.e.~$\theta$,
\begin{equation} \lb{4.15}
\Real f(e^{i\theta}) =g(e^{i\theta})
\end{equation}
\end{theorem}

\begin{proof} Given $\theta_0\in (0,\pi)$, break the integral in \eqref{4.14} into
two parts: $I_1\equiv (\theta_0 -\delta, \theta_0+\delta)\cup (-\theta_0-\delta,
-\theta_0 +\delta)$ for $\abs{\delta} <\min (\theta_0, \pi-\theta_0)$ and the
complement, $I_2$. By \eqref{4.12}, if $\theta\in I_2$,
\begin{equation} \lb{4.16}
\abs{Q(re^{i\theta_0}, e^{i\theta}}\leq C\sin^2\theta
\end{equation}
uniformly in $r$ and $\lim_{r\uparrow 1} Q(re^{i\theta_0}, e^{i\theta})$ exists
and is pure imaginary. Thus the part of the integral in \eqref{4.14} for
$\theta\in I_2$ has a limit with real part $0$; if $z=re^{i\theta_0}$,
$r\uparrow 1$.

On $I_1$, we can rewrite $Q$ as a sum of its four summands $(\f12 P(z,e^{i\theta}),
\f12 P(z, e^{-i\theta}), \alpha(z)$, and $\beta(z)\cos\theta$). Clearly, $\alpha
(re^{i\theta_0})$ and $\beta (re^{i\theta_0})$ have limits which are pure imaginary.
By the standard theory of Poisson kernels (Rudin \cite{Rudin}, Duren \cite{Duren}),
the $P$ terms have a limit for a.e.~$\theta_0$ whose real part is $\f12(g(e^{i\theta_0})
+ g(e^{-i\theta_0})) =g(e^{i\theta_0})$ by the assumed symmetry of $g$.
\end{proof}

\section{A Necessary and Sufficient Condition\\for Jost Asymptotics} \lb{s5}

Our goal in this section is to prove that

\begin{theorem}\lb{T5.1} Let $J$ be a Jacobi matrix with $a_n\to 1$, $b_n\to 0$.
Let $Q=\{z\in\bbD\mid z+z^{-1}\text{ is an eigenvalue of $J$}\}$. Then the
following are equivalent:
\begin{SL}
\item[{\rm{(i)}}] Szeg\H{o} asymptotics {\rm{(}}i.e., $z^n p_n (z+\f{1}{z})$
converges to a nonzero limit as $n\to\infty${\rm{)}} hold for all $z\in\bbD
\backslash Q$ uniformly on compact subsets of $\bbD\backslash Q$.

\item[{\rm{(ii)}}] Szeg\H{o} asymptotics hold for all $z$ with $\abs{z}=\veps$
for some $\veps >0$ and uniformly in such $z$.

\item[{\rm{(iii)}}] Jost asymptotics {\rm{(}}i.e., $z^{-n} w_n(z)$ has a
nonzero limit{\rm{)}} hold for all $z\in\bbD\backslash Q$ uniformly on
compact subsets of $\bbD\backslash Q$.

\item[{\rm{(iv)}}] Jost asymptotics hold for all $z$ with $\abs{z}=\veps$
for some $\veps >0$ uniformly in such $z$.

\item[{\rm{(v)}}] The $a$'s and $b$'s obey three conditions:
\begin{SL}
\item[$(\alpha)$]
\begin{equation}
\sum_{n=1}^\infty \, \abs{a_n-1}^2 + \sum_{n=1}^\infty \, \abs{b_n}^2 <\infty \lb{5.1}
\end{equation}
\item[$(\beta)$] $\lim_{n\to\infty} \, a_1 \dots a_n$ exists and is not zero.
\item[$(\gamma)$] $\lim_{n\to\infty} \, \sum_{j=1}^n b_j$ exists.
\end{SL}
\item[{\rm{(vi)}}] The spectral measure, $\mu$, on $\bbR$ and orthonormal
polynomials obey the following properties:
\begin{SL}
\item[$(\delta)$] $\int_{-2}^2 \log ( {d\mu_{\ac}/dE}) \sqrt{4-E^2} \, dE >-\infty$
\item[$(\veps)$] $\sum_n ( \abs{E_n^\pm} -2)^{3/2} <\infty$
\item[$(\kappa)$] If the orthonormal polynomials have the form
\begin{equation} \lb{5.1a}
p_n(x) =\gamma_n (x^n -\lambda_n x^{n-1} + \cdots )
\end{equation}
then $\lim_{n\to\infty} \gamma_n$ exists and is nonzero and $\lim_{n\to\infty}
\lambda_n$ exists.
\end{SL}
\end{SL}
\end{theorem}

\begin{remarks} 1. We will see shortly that $w_n(z)$ has an $n$th-order zero at $z=0$,
so $z^{-n} w_n(z)$ has a removable singularity at $z=0$ --- and it is that value
we intend when we say the limit exists at $z=0$.

\smallskip
2. We will discuss below what happens at the $z_0$'s in $Q$. (Basically, $z^n p_n
(z+\f{1}{z}$) has a zero limit there and, by shifting from Weyl to Jost solutions,
we will also have control at $z_0$'s in $Q$ of the other solutions.)

\smallskip
3. We will see that $u(z)\equiv (\lim w_\infty (z))^{-1}$ always has a
factorization formula when (v) holds. $u$ will be expressed in terms of
``spectral data" and the limits in $(\beta)$ and $(\gamma)$.
\end{remarks}

Define
\begin{equation} \lb{5.2}
M(z,J)=\langle \delta_1, (z+z^{-1} -J)^{-1} \delta_1 \rangle =w_1(z)
\end{equation}
Let $J^{(n)}$ be the Jacobi matrix obtained by removing the first $n$ rows and left $n$
columns of $J$. Let
\begin{equation} \lb{5.3}
M_n(z,J) =M(z, J^{(n)})
\end{equation}
so $M_0 (z,J) \equiv M(z,J)$. We will often drop the $J$ if it is fixed in some discussion.

\begin{lemma} \lb{L5.2}
\begin{alignat}{2}
&\text{\rm{(i)}} \qquad && M(z) =z+O(z^2) \lb{5.4} \\
&\text{\rm{(ii)}} \qquad && M_n(z) = \f{w_{n+1}(z)}{a_n w_n (z)} \lb{5.5} \\
&\text{\rm{(iii)}} \qquad && w_n(z) =M(z) (a_1 M_1 (z)) \dots (a_{n-1} M_{n-1} (z)) \lb{5.6} \\
&\text{\rm{(iv)}} \qquad && w_n(z) =(a_1 \dots a_n) z^n + O(z^{n+1}) \lb{5.7} \\
&\text{\rm{(v)}} \qquad && M_n(z) =(z+z^{-1} -b_{n+1} -a_{n+1}^2 M_{n+1}(z))^{-1} \lb{5.8} \\
&\text{\rm{(vi)}} \qquad && \log \biggl(\f{M_n (z)}{z}\biggr) = b_{n+1} z + (\tfrac12\, b_{n+1}^2 +
a_{n+1}^2 -1) z^2 + O(z^3) \lb{5.9}
\end{alignat}
\end{lemma}

\begin{remark} Some of these equalities are intended in the sense of the
field of meromorphic functions. For example, if $\ell <n$ and $w_\ell (z_0)=0$, then
$M_\ell (z)$ has a pole at $z_0$ and $M_{\ell-1}(z)$ a zero there and they are intended
to cancel in \eqref{5.6}. Alternatively, these formulae hold initially away from
$\{z\in\bbD\mid z+z^{-1} \in\sigma (J^{(\ell)})$ for some $\ell=0,1,\dots\}$ and then
they have removable singularities in some cases.
\end{remark}

\begin{proof}
\begin{SL}
\item[(i)] $M(z)/z =\langle \delta_1, (1+z^2 -zJ)^{-1}\delta_1\rangle = 1 +O(z)$ as $z\to 0$.

\item[(ii)] As noted in Section~\ref{s2}, $w_n(z)$ is normalized by \eqref{2.5a}, that is, by
\[
a_1 w_2 (z) + (b_1 -z-z^{-1}) w_1 (z) =-1
\]
and, of course, $M(z) =w_1(z)$.  \eqref{5.5} thus follows from
\[
a_{n+1} \biggl( \f{w_{n+2}}{a_n w_n} \biggr) + (b_{n+1} -z-z^{-1})
\biggl( \f{w_{n+1}}{a_n w_n}\biggr)
= -1
\]
since $w_{n+j}/a_n w_n$ solves the difference equation for $J^{(n)}$.

\item[(iii)] follows from \eqref{5.5} and $w_1 =M$.

\item[(iv)] is immediate from \eqref{5.4} for $M_n(z)$ and \eqref{5.6}.

\item[(v)] follows from \eqref{5.5} and the difference equation for $w$.

\item[(vi)] From \eqref{5.8} for $n=0$ and \eqref{5.4} for $M_1$,
\[
\f{M(z)}{z} = (1-b_1 z -a_1 z^2 + z^2 + O(z^3))^{-1}
\]
so
\begin{align*}
\log \biggl( \f{M(z)}{z}\biggr) &= -\log (1-b_1 z - a_1^2 z^2 + z^2 +O(z^3)) \\
&= (b_1 z + \tfrac12\, b_1 z^2) + (a_1^2 -1) z^2 + O(z^3)
\end{align*}
\end{SL}
\end{proof}

\begin{proof}[Reduction of Theorem~\ref{T5.1} to {\rm{(v)}} $\Rightarrow$ {\rm{(iii)}}]
By Theorem~\ref{T2.2}, (i) $\Leftrightarrow$ (iii) and (ii) $\Leftrightarrow$ (iv).
(iii) $\Rightarrow$ (iv) is trivial. Thus we need to prove (iv) $\Rightarrow$ (v)
and (v) $\Leftrightarrow$ (vi) to reduce the proof to (v) $\Rightarrow$ (iii).

The equivalence of (v) and (vi) is easy, given the result of Killip-Simon \cite{KS}.
They prove that $(\alpha)$ $\Leftrightarrow$ $(\delta), (\veps)$. The equivalence of
$(\kappa)$ and $(\beta),(\gamma)$ is immediate since the recursion relations for $p$
imply that
\begin{align*}
a_{n+1} \gamma_{n+1} &= \gamma_n \\
-\lambda_{n+1} &= -\lambda_n - b_{n+1}
\end{align*}
so
\begin{equation} \lb{5.9a}
\gamma_n =(a_1 \dots a_n)^{-1} \qquad \lambda_n =\sum_{j=1}^n b_j
\end{equation}

To study (iv) $\Rightarrow$ (v), define
\begin{equation} \lb{5.10}
\ti w_n(z) =z^{-n} w_n(z) \qquad \ti M_n(z) =z^{-1} M_n(z)
\end{equation}
so, by \eqref{5.6},
\begin{equation} \lb{5.11}
\log \ti w_n(z) =\sum_{j=1}^{n-1} \log (a_j) + \sum_{j=1}^{n-1} \log \ti M_{j-1}(z)
\end{equation}
Convergence of $\ti w_n(z)$ uniformly on the circle and analyticity of $\ti w_n(z)$
implies the derivatives of $\ti w_n(z)$ at $z=0$ all converge. By \eqref{5.9} and
\eqref{5.11}, the terms of order $1,z,z^2$ yield
\begin{gather}
\lim_{n\to\infty} \, \sum_{j=1}^{n-1} \log (a_j) =\nu_1 \lb{5.12} \\
\lim_{n\to\infty} \, \sum_{j=1}^{n-1} b_j = \nu_2 \lb{5.13} \\
\lim_{n\to\infty} \, \sum_{j=1}^{n-1} a_j^2 - 1+\tfrac12\, b_j^2 = \nu_3 \lb{5.14}
\end{gather}
all exist. Following Killip-Simon \cite{KS}, we look at $\eqref{5.14} -2\times
\eqref{5.12}$ to see
\begin{equation} \lb{5.15}
\lim_{n\to\infty} \, \sum_{j=1}^{n-1} G(a_j) + \tfrac12\, b_j^2 =\nu_3 -2\nu_1
\end{equation}
where
\[
G(a) =a^2 -1-2\log (a)
\]
Since $G(a) >0$ for $a\in (0,\infty)$, the summand in \eqref{5.15} is nonnegative,
so \eqref{5.15} implies absolute convergence. Since $G(a)\geq (a-1)^2$ (for $G(1)
=G'(1)=0$ and $G''(a)\geq 2$), \eqref{5.15} implies
\[
\sum_{j=1}^\infty (a_j-1)^2 + \tfrac12\, b_j^2 <\infty
\]
which is $(\alpha)$. $(\beta)$ is the exponential of \eqref{5.12} and $(\gamma)$
is \eqref{5.13}.
\end{proof}

We now turn towards proving Jost asymptotics when $(\alpha)$, $(\beta)$, $(\gamma)$
hold. We will give three proofs: one in this section using canonical factorization
of $M$-functions, one in Section~\ref{s6} using renormalized determinants, and one
in Section~\ref{s7} using Levinson-type asymptotic analysis of difference equations.

Our starting point for the proof in this section will be the ``nonlocal"
step-by-step sum rule of Simon \cite{Sim288}:

\begin{theorem} \lb{T5.3} For any Jacobi matrix with $a_n\to 1$, $b_n\to 0$,
\begin{equation} \lb{5.16}
\{\theta\mid\Ima M_n(\theta) \neq 0\} =\{\theta\mid\Ima M_{n+1}(\theta)\neq 0\}
\end{equation}
{\rm{(}}modulo sets of $d\theta/2\pi$-measure zero{\rm{)}}. For any $p<\infty$,
\begin{equation} \lb{5.17}
\log\biggl( \f{\Ima M_n (\theta)}{\Ima M_{n-1}(\theta)}\biggr) \in L^p
\biggl(\partial\bbD, \f{\partial\theta}{2\pi}\biggr)
\end{equation}
and
\begin{equation} \lb{5.18}
a_{n+1} M_n(z) = zB_n^+ B_n^- (z) \exp \biggl( \f{1}{4\pi} \int \log
\biggl( \f{\Ima M_n(z)}{\Ima M_{n+1}(z)} \biggr)
\biggl(\f{e^{i\theta}+z}{e^{i\theta}-z}\biggr) d\theta\biggr)
\end{equation}
Here $B_n^\pm$ are alternating Blaschke products {\rm{(}}$B^+$ for $0<p_{1,+}^{(n)}
<z_{1,+}^{(n)} < \cdots < p_{\ell,+}^{(n)} < z_{\ell,+}^{(n)}$ and $B^-$ for
$0>p_{1,-}^{(n)} > z_{1,-}^{(n)} > \cdots${\rm{)}} with $p_{j,\pm}^{(n)} +
p_{j,\pm}^{(n)-1}$ the eigenvalues of $J^{(n)}$ and $z_{j,\pm}^{(n)} +
z_{j,\pm}^{(n)-1}$ the eigenvalues of $J^{(n+1)}$.
\end{theorem}

\begin{remarks} 1. \eqref{5.18} is a special case of a general factorization
theorem for meromorphic Herglotz functions, $f$, of $\bbD$. The general theorem
has $\f{1}{2\pi} \log (\abs{f(e^{i\theta})})$. \eqref{5.18} then follows from
$\abs{a_{n+1} M_n}^2 =\Ima M_n/\Ima M_{n+1}$, which is a consequence of \eqref{5.8}.

\smallskip
2. In our applications, the set in \eqref{5.16} is all of $\partial\bbD$. When
this is false, $\Ima M_n/ \Ima M_{n+1}$ in \eqref{5.17} and \eqref{5.18} have
to be suitably defined on the complement of the set in \eqref{5.16}; see
\cite{Sim288} for details.
\end{remarks}

\smallskip
We define
\begin{align}
L_n(z) &=\log \biggl( \f{a_{n+1} M_n(z)}{z}\biggr) \lb{5.19}  \\
N_n(z) &= L_n (z) -\alpha (z) L_n(0) -\tfrac12\, \beta(z) L'_n(0) \lb{5.20}
\end{align}
where $\alpha,\beta$ are given by \eqref{3.16}. If $p_{1,\pm}^{(n)}$ are the
poles of $M_n$ closest to $z=0$, we define $L_n(z)$ unambiguously on $\bbD
\backslash [p_{1,+}^{(n)}, 1) \cup (-1, p_{1,-}^{(n)}]$ by requiring $L_n(z)$
analytic and $L_n(0)$ real. Since $p_{1,\pm}^{(n)}\to \pm 1$ as $n\to\infty$,
the result below exponentiated holds on $\bbD\backslash \{p_{j,\pm}^{(0)}
\mid j=1,2,\dots\}$.

\begin{lemma}\lb{L5.4} Suppose that $(\alpha)$, that is, \eqref{5.1} holds.
Then for all $z\in\bbD \backslash [p_{1,+}^{(0)}, 1)\cup (-1, p_{1,+}^{(0)}] =S$,
\begin{equation} \lb{5.21}
\lim_{N\to\infty} \, \sum_{n=0}^N N_n(z)
\end{equation}
exists and the convergence is uniform on compact subsets of $S$.
\end{lemma}

\begin{proof} By \eqref{5.18},
\[
L_n(z) =\log B_n^+ (z) + \log B_n^-(z) + \f{1}{4\pi} \int
\biggl( \f{e^{i\theta}+z}{e^{i\theta}-z}\biggr)
\log \biggl( \f{\Ima M_n}{\Ima M_{n+1}}\biggr)\, d\theta
\]
Using  \eqref{3.19} and \eqref{4.4},
\begin{equation} \lb{5.23}
\begin{split}
N_n(z) &=\sum_{j,\pm} [-\log q(z,p_{j,\pm}^{(n)})+\log q (z,z_{j,\pm}^{(n)})] \\
&\qquad + \f{1}{4\pi} \int Q(z, e^{i\theta}) \log \biggl( \f{\Ima M_n}{\Ima
_{n+1}}\biggr)\, d\theta
\end{split}
\end{equation}

Since $z_{j,\pm}^{(n)}=p_{j,\pm}^{(n+1)}$, this implies
\begin{equation} \lb{5.24}
\begin{split}
\sum_{n=0}^N N_n(z) &=\sum_{j,\pm} [-\log q(z, p_{j,\pm}^{(0)}) +
\log q(z, p_{j,\pm}^{(N+1)})] \\
&\qquad + \f{1}{4\pi} \int Q(z, e^{i\theta}) \log \biggl( \f{\Ima M}{\Ima M_{N+1}}\biggr)\,
d\theta
\end{split}
\end{equation}

The Killip-Simon \cite{KS} $P_2$ sum rule implies
\begin{SL}
\item[(a)] $\int \sin^2\theta \abs{\log (\f{\Ima M}{\sin\theta})} \f{d\theta}{2\pi}<\infty$
\item[(b)] $\sum (1-\abs{p_{j,\pm}^{(0)}})^3 <\infty$
\item[(c)] $\lim_{N\to\infty} \int \sin^2\theta \abs{\log (\f{\Ima M_{N+1}}{\sin\theta})}
\f{d\theta}{2\pi} =0$
\item[(d)] $\lim_{N\to\infty} \sum (1-\abs{p_{j,\pm}^{(N+1)}})^3 =0$
\end{SL}

(a) and (b) and the estimates \eqref{3.23} and \eqref{4.13} allow us to write
\eqref{5.24} as a difference of $M^{(0)}/p^{(0)}$ terms and $M^{(N+1)}/p^{(N+1)}$
terms and then (c),(d) show that the error terms go to zero. The result is that
$\lim \sum_{n=0}^N N_n(z)$ exists and
\begin{equation} \lb{5.25}
\begin{split}
\lim_{n\to\infty} &\exp \biggl(\, \sum_{n=0}^N N_n(z)\biggr) \\
& = \prod_{j=\pm 1} q(z,p_{j,\pm}^{(0)})^{-1}
\exp\biggl( \f{1}{4\pi}\int Q(z, e^{i\theta}) \log
\biggl( \f{\Ima M}{\sin\theta}\biggr)\, d\theta\biggr)
\end{split}
\end{equation}
The proof shows the convergence is uniform.
\end{proof}

\begin{proof}[Proof of Theorem~\ref{T5.1}{\rm{(v)}} $\Rightarrow$ {\rm{(iii)}}]
By \eqref{5.6}, with $\ti w_n(z) = z^{-n} w_n(z)$,
\begin{equation} \lb{5.26}
a_n \ti w_n(z) =\exp \biggl(\, \sum_{j=0}^{n-1} L_j(z)\biggr)
\end{equation}
Since $a_n\to 1$, $\ti w_n(z)$ has a nonzero limit (i.e., Jost asymptotics hold)
if and only if
\[
\lim_{N\to\infty} \, \sum_{j=0}^N L_j(z)
\]
exists and the convergence of $\ti w_n$ is uniform if and only if the convergence
of the sum is. Since
\[
L_j(0) =\log (a_{j+1})
\]
by \eqref{5.4} and
\[
L'_j(0) =b_{j+1}
\]
by \eqref{5.9}, we have that
\begin{equation} \lb{5.27x}
L_j(z) =N_j(z) +\alpha(z) \log (a_{j+1}) + \tfrac12\, \beta(z) b_{j+1}
\end{equation}

By Lemma~\ref{5.3}, $\sum_{n=0}^N N_n(z)$ converges uniformly if $(\alpha)$ holds,
$(\beta)$ and $(\gamma)$ say that $\sum_{n=0}^N \log (a_{n+1})$ and $\sum_{n=0}^N
b_{j+1}$ converge so $\sum_{n=0}^N L_n(z)$ converges uniformly.
\end{proof}

If $\ti w_n(z)$ has a nonzero limit, we define the Jost function by
\begin{equation} \lb{5.27}
u(z) =\lim_{n\to\infty} \, \ti w_n(z)^{-1}
\end{equation}
This agrees with the usual definition if $\sum_{n=1}^\infty \abs{a_n-1} +
\abs{b_n}<\infty$. Thus
\begin{equation} \lb{5.28}
u(z) =\exp\biggl( -\sum_{n=0}^\infty L_n(z)\biggr)
\end{equation}
and we have proven (by \eqref{5.25}) that

\begin{theorem}\lb{T5.5} If $(\alpha)$, $(\beta)$, $(\gamma)$ hold, then
\begin{equation} \lb{5.29}
\begin{split}
u(z) & =\biggl(\, \prod_{j=1}^\infty a_j\biggr)^{-\alpha(z)} e^{-\f12 \beta(z)
\sum_{j=1}^\infty b_j} \\
&\qquad \prod_{j=1,\pm}^\infty q(z,p_{j,\pm}^{(0)}) \exp\biggl( -\f{1}{4\pi}
\int Q(z,e^{i\theta})\log \biggl( \f{\Ima M}{\sin\theta}\biggr) \, d\theta\biggr)
\end{split}
\end{equation}
\end{theorem}

In the above, $\prod_{j=1}^\infty a_j$ and $\sum_{j=1}^\infty b_j$ refer to the
conditional limits.

The integral representation \eqref{5.29} implies

\begin{theorem}\lb{T5.6} Let $(\alpha)$, $(\beta)$, $(\gamma)$ hold and let
\[
u(z) =\bigl(\, \lim_{n\to\infty}\, z^{-n} w_n(n)\bigr)^{-1}
\]
Then
\begin{SL}
\item[{\rm{(i)}}] After removing the removable singularities at $\{p_j^{(0)}\}$, $u(z)$
is analytic in $\bbD$ and $u(z_0)=0$ {\rm{(}}$z_0\in\bbD${\rm{)}} if and only if
$z_0\in\{p_j^{(0)}\}$, that is, if and only if $z_0 + z_0^{-1}$ is an eigenvalue of $J$.

\item[{\rm{(ii)}}] For a.e.~$\theta$, $\lim_{r\uparrow 1} u(re^{i\theta})\equiv
u(e^{i\theta})$ exists and
\begin{equation} \lb{5.30}
\Ima M(e^{i\theta}) \abs{u(re^{i\theta})}^2 =\sin\theta
\end{equation}
for a.e.~$\theta$.
\end{SL}
\end{theorem}

\begin{proof} (i) is immediate from \eqref{5.29} and (ii) follows from the fact
that $\alpha,\beta$ have purely imaginary values on $\partial\bbD\backslash\{\pm 1\}$
from \eqref{3.30} and Theorem~\ref{T4.2}.
\end{proof}

\begin{remark} However, unlike the case where $\sum_{n=1}^\infty \abs{a_n-1}
+\abs{b_n}<\infty$, $u$ may not be Nevanlinna. Indeed, if $\sum_{n=1}^\infty
(\abs{E_n^\pm}-2)^{1/2}=\infty$, $u$ cannot be Nevanlinna. This is the subject
of Section~\ref{s8}.
\end{remark}

\section{Renormalized Determinants} \lb{s6}

The idea that Jost functions are given by Fredholm determinants goes back to
Jost-Pais \cite{JP}, and for Jacobi matrices was made explicit by Killip-Simon \cite{KS}.
They define the perturbation determinant by
\begin{equation} \lb{6.1}
L(z,J)=\det (\boldsymbol{1} + \delta J(J_0 -(z+z^{-1}))^{-1})
\end{equation}
where
\begin{equation} \lb{6.2}
\delta J =J-J_0
\end{equation}
and $J_0$ is the Jacobi matrix associated to $a_n \equiv 1$, $b_n\equiv 0$. This definition
is used when $z\in\bbD$ and
\begin{equation} \lb{6.3}
\sum_{n=1}^\infty\, \abs{a_n -1} + \abs{b_n} < \infty
\end{equation}
In this case, $\delta J$ is trace class and the $\det$ in \eqref{6.1} is the standard
trace class determinant (see Simon \cite{TI} and Goh'berg-Krein \cite{GK}).

What Killip-Simon \cite{KS} prove in their Theorem~2.16 is

\begin{theorem} \lb{T6.1} For $J-J_0$ trace class, $z\in\bbD$, and $z+z^{-1} \notin
\sigma(J)$, we have that the function $M$\! given by \eqref{5.3} obeys
\begin{equation} \lb{6.4}
M(z,J) = \f{zL(z,J^{(1)})}{L(z,J)}
\end{equation}
and with $w_n(z)$, the Weyl solution,
\begin{equation} \lb{6.5}
z^{-n} w_n(z) \to \f{[\prod_{j=1}^\infty a_j]}{L(z,j)}
\end{equation}
\end{theorem}

\begin{remarks} 1. \eqref{6.4} implies \eqref{6.5} using \eqref{5.6} and $L(z,J^{(n)})
\to 1$ since $\|J^{(n)} -J_0\|_1 \to 0$.

\smallskip
2. \eqref{6.5}, of course, says that the Jost function is given by
\begin{equation} \lb{6.6}
u(z;J) = \biggl[\, \prod_{j=1}^\infty a_j \biggr]^{-1} L(z;J)
\end{equation}

\smallskip
3. Formula \eqref{6.4} is an expression of Cramer's rule since, very formally, Cramer's
rule says
\begin{equation} \lb{6.7}
M(z,J) = \f{\det(z+z^{-1} -J^{(1)})}{\det (z+z^{-1} -J)}
\end{equation}
and
\begin{equation} \lb{6.8}
z=M(z,J^{(0)}) = \f{\det (z+z^{-1} -J_0^{(1)})}{\det (z+z^{-1} -J_0)}
\end{equation}
Moreover,
\begin{equation} \lb{6.9}
L(z,J) = \f{\det (J-(z+z^{-1}))}{\det (J_0 - (z+z^{-1}))}
\end{equation}
\end{remarks}

\eqref{6.7}--\eqref{6.9} manipulate to \eqref{6.4}. Of course, the $\det$'s in
\eqref{6.7}--\eqref{6.9} are all infinite, but one way to prove \eqref{6.4} is to
prove \eqref{6.7}--\eqref{6.9} for cutoff finite matrices and take limits.

When \eqref{5.1} holds, $J-J_0$ may not any longer be trace class, but it is
Hilbert-Schmidt, which suggests that we use the renormalized determinant for such
operators. Such determinants go back to Carleman \cite{ti25}. They are discussed
in \cite{TI,GK}. Our approach, due to Seiler \cite{ti142,ti144} and used in \cite{TI},
relies on the fact that if $A$ is Hilbert-Schmidt, then
\[
B=(1+A)e^{-A} -\boldsymbol{1}
\]
is trace class, so we can define
\begin{equation} \lb{6.10}
\det_2 (\boldsymbol{1} +A) \equiv \det (\boldsymbol{1} +B) =
\det ((\boldsymbol{1}+A) e^{-A})
\end{equation}

It obeys (see \cite[Chapter~9]{TI})
\begin{gather}
A\in\calI_1 \rightarrow \det (\boldsymbol{1}+A) = \det_2 (\boldsymbol{1}+A)
e^{\tr (A)} \lb{6.11} \\
\abs{\det_2 (\boldsymbol{1} + A) -\det_2 (\boldsymbol{1}+C)} \leq
\|A-C\|_2 \exp (\Gamma_2 (\|A\|_2 + \|C\|_2 + 1))  \lb{6.12}
\end{gather}
for a suitable constant, $\Gamma_2$.

We note (see \cite[eqn. (1.2.24)]{OPUC1}) that
\begin{equation} \lb{6.13}
[(J_0 - (z+z^{-1}))^{-1}]_{nm} = -(z^{-1}-z)^{-1} [z^{\abs{m-n}} -z^{m+n}]
\end{equation}
Thus
\begin{equation} \lb{6.14}
[\delta J (J_0 -(z+z^{-1}))^{-1}]_{nn} = -(z^{-1} -z)^{-1} \{(a_{n-1} -1) (z-z^{2n-1}) +
b_n (1-z^{2n}) + (a_n -1) (z-z^{2n+1})\}
\end{equation}
which implies:

\begin{lemma} \lb{L6.2} If $\delta J$ is trace class and $z\in\bbD$, then
\begin{equation} \lb{6.15}
\tr (\delta J (J_0 -(z+z^{-1}))^{-1}) = -(z^{-1} -z)^{-1}
\biggl\{\, \sum_{n=1}^\infty b_n (1-z^{2n}) + 2\sum_{n=1}^\infty (a_n -1)
(z-z^{2n+1})\biggr\}
\end{equation}
\end{lemma}

It also explains the relevance of

\begin{proposition} \lb{P6.3} Suppose $a_n, b_n$ obey $(\alpha)$--$(\gamma)$ of
Theorem~\ref{T5.1}. Then
\begin{equation} \lb{6.16}
T(z;J) = \lim_{N\to\infty}\, \biggl[ -(z^{-1} -z)^{-1} \sum_{n=1}^N b_n (1-z^{2n}) +
2\sum_{n=1}^N (a_n -1) (z-z^{2n+1})\biggr]
\end{equation}
exists for all $z\in\bbD$ and the convergence is uniform for compact subsets of $\bbD$.
\end{proposition}

\begin{proof} $z^{2\boldsymbol{\cdot}}\in\ell^2$ for $z\in\bbD$ uniformly on compact subsets, so
$(\alpha)$ implies
\[
\sum_{n=1}^\infty b_n z^{2n} \qquad\text{and} \qquad \sum_{n=1}^\infty (a_n -1) z^{2n+1}
\]
converge absolutely to an analytic limit.

$(\beta)$ plus $\sum_{n=1}^\infty \abs{a_n -1}^2 <\infty$ implies $\lim_{N\to\infty}
\sum_{n=1}^N (a_n -1)$ exists, and this plus $(\gamma)$ implies the existence of the
remaining terms.
\end{proof}

\begin{definition} If $(\alpha)$--$(\gamma)$ of Theorem~\ref{T5.1} hold, we
define
\begin{equation} \lb{6.17}
L_\ren (z,J) = \det_2 (1+\delta J(J_0 - (z+z^{-1}))^{-1}) e^{T(z;J)}
\end{equation}
\end{definition}

\begin{proposition}\lb{P6.3A} Let $z\in\bbD$. $L_\ren (z,J)=0$ if and only if $z+z^{-1}
\in \sigma(J)$.
\end{proposition}

\begin{proof} $\det_2 (1+A)=0$ if and only if $1+A$ is not invertible (see
\cite[Chapter~9]{TI}). Since $1+\delta J(J_0 - (z+z^{-1}))^{-1} = (J-(z+z^{-1}))
(J_0 - (z+z^{-1}))$, this happens if and only if $z+z^{-1}\in \sigma (J)$.
\end{proof}

By \eqref{6.1}, \eqref{6.11}, and \eqref{6.15}, we have that

\begin{proposition} \lb{P6.4} If $\sum_{n=1}^\infty \abs{a_n -1} + \abs{b_n}<\infty$,
then
\begin{equation} \lb{6.18}
L_\ren (z,J) = L(z,T)
\end{equation}
\end{proposition}

\begin{theorem} \lb{T6.5} If $(\alpha)$--$(\gamma)$ of Theorem~\ref{T5.1} hold, then
for all $z\in\bbD$ so $z+z^{-1}\notin\sigma (J)$,
\begin{equation} \lb{6.19}
M(z,J) = \f{z L_\ren (z;J^{(1)})}{L_\ren (z;J)}
\end{equation}
and
\begin{equation} \lb{6.20}
z^{-n} w_n (z) = \biggl( \, \prod_{j=1}^{n-1} a_j\biggr) \,
\f{L_\ren (z;J^{(n)})}{L_\ren (z;J)}
\end{equation}
\end{theorem}

\begin{proof} \eqref{6.20} is implied by \eqref{6.19} and \eqref{5.6}, so we need only
prove \eqref{6.19}. Define $J_n$ to be the Jacobi matrix with
\begin{align}
a_j^{(n)} &= \begin{cases}
a_j & j\leq n-1 \\
1 & j \geq n
\end{cases} \lb{6.21} \\
b_j^{(n)} &= \begin{cases}
b_j & j\leq n \\
0 & j\geq n+1
\end{cases} \lb{6.22}
\end{align}
Then, by $(\alpha)$,
\begin{equation} \lb{6.23}
\|J_n -J\|_2 \to 0
\end{equation}
so $\|\delta J_n -\delta J\|_2 \to 0$, so by \eqref{6.12},
\begin{equation} \lb{6.24}
\det_2 (1+\delta J_n (J_0 - (z+z^{-1}))^{-1}) \to \det_2
(1+\delta J (J_0 - (z+z^{-1}))^{-1})
\end{equation}
It is easy to see that $T(z,J_n) \to T(z,J)$. Thus, using \eqref{6.18}, uniformly on
compacts of $\bbD$,
\begin{equation} \lb{6.25}
\lim_{n\to\infty} \, L(z,J_n) = \lim_{n\to\infty} \,
L_\ren (z,J_n) = L_\ren (z,J)
\end{equation}
The same is true for $J_n^{(1)}$ and $J^{(1)}$. Therefore, since
$u(z,J)\neq 0$, \eqref{6.4} implies \eqref{6.19}.
\end{proof}

We therefore have the second proof of the hard part of Theorem~\ref{T5.1}:

\begin{theorem} \lb{T6.6} Let $(\alpha)$--$(\gamma)$ of Theorem~\ref{T5.1}
hold. Then uniformly on compact subsets of $\bbD\backslash \{z\mid z+z^{-1}
\in \sigma(J)\}$,
\begin{equation} \lb{6.26}
z^{-n} w_n(z) \to u(z)^{-1}
\end{equation}
where
\begin{equation} \lb{6.27}
u(z) = \biggl(\, \lim_{n\to\infty} \, \prod_{j=1}^n a_j \biggr)^{-1}
L_\ren (z,J)
\end{equation}
\end{theorem}

\begin{proof} By \eqref{6.20}, this is equivalent to
\begin{equation} \lb{6.28}
\lim_{n\to\infty} \, L_\ren (z,J^{(n)}) =1
\end{equation}
uniformly on compacts. It is easy to see that $\lim_{n\to\infty} T(z,J^{(n)})=0$.
Since $\|J^{(n)} -J_0\|_2 \to 0$,
\[
\det_2 (1+\delta J_n (J_0 -(z+z^{-1}))^{-1}) \to 1
\]
This proves \eqref{6.28}.
\end{proof}

\section{Geronimo-Case Equations} \lb{s7}

Given a set $\{a_n, b_n\}_{n=1}^\infty$ of real Jacobi parameters, the
Geronimo-Case polynomials $c_n(z), g_n(z)$ are defined by the recursion
relations:
\begin{align}
c_{n+1}(z) &= a_{n+1}^{-1} [(z^2 -b_{n+1} z) c_n(z) + g_n(z)] \lb{7.1} \\
g_{n+1}(z) &= a_{n+1}^{-1} [((1-a_{n+1}^2)z^2 -b_{n+1}z) c_n(z) + g_n(z)] \lb{7.2}
\end{align}
with initial conditions
\begin{equation} \lb{7.3}
c_0 (z) = g_0(z) =1
\end{equation}
They were introduced in a slightly different form by Geronimo-Case \cite{GC80}
who, under a condition that $\sum n[\abs{a_n -1} + \abs{b_n}] <\infty$, proved
that for $z\in\ol{\bbD}$, $\lim_{n\to\infty} g_n(z)$ exists and defined it to
be the Jost function. In Paper II of our current series \cite{Jost2}, we will reexamine
these equations to prove convergence in $\bbD$ if $\sum_n \abs{a_n-1} + \abs{b_n}
<\infty$ and, most importantly, identify what $c_n$ and $g_n$ are, namely,
\begin{equation} \lb{7.4}
c_n(z) = z^n p_n \biggl( z+\f{1}{z}\biggr)
\end{equation}
where $p_n$ are the orthonormal polynomials. Moreover, if $\ti J_\ell$ is defined
like $J_\ell$ (see \eqref{6.21}/\eqref{6.22}) but with a different cutoff on
$a_j$, that is,
\begin{align}
\ti a_j^{(\ell)} &= \begin{cases}
a_j & j\leq \ell \\
1 & j \geq \ell+1
\end{cases} \lb{7.5}  \\
\ti b_j^{(\ell)} &= \begin{cases}
b_j & j\leq \ell \\
0 & j\geq \ell+1
\end{cases} \lb{7.6}
\end{align}
then
\begin{equation} \lb{7.7}
g_n (z, J) = u (z,\ti J_n)
\end{equation}
the conventional Jost function for $\ti J_n$, that is, $(\lim_{m\to\infty} z^{-m}
w_m (z,\ti J_n))^{-1}$.

Our goal here is to extend Theorem~\ref{T5.1} by proving:

\begin{theorem} \lb{T7.1} The following are equivalent:
\begin{SL}
\item[{\rm{(a)}}] For some $\veps\in (0,1)$, $\lim_{n\to\infty} c_n (z)$ exists for
$\abs{z}=\veps$ uniformly in such $z$.

\item[{\rm{(b)}}] For all $z\in\bbD$, $\lim_{n\to\infty} c_n(z)$ and $\lim_{n\to\infty}
g_n(z)$ exist uniformly on compacts of $\bbD$, and $\lim g_n(z)$ is the Jost function
$u(z;J)$.

\item[{\rm{(c)}}] Conditions $(\alpha)$--$(\gamma)$ of Theorem~\ref{T5.1} hold.
\end{SL}
\end{theorem}

\begin{proof} That (a) $\Rightarrow$ (c) is just (ii) $\Rightarrow$ (iv)
in Theorem~\ref{T5.1}, and (b) $\Rightarrow$ (a) is trivial. So we only need
(c) $\Rightarrow$ (b). Convergence of $c_n$ is just (iv) $\Rightarrow$ (i)
of Theorem~\ref{T5.1}, so we only need convergence of $g_n$. To see this, we
use \eqref{7.7}, \eqref{6.17}, and \eqref{6.27}. $\|\ti J_n -J\|_2 \to 0$, so
\begin{equation} \lb{7.8}
\det_2 (1+\delta\ti J_n (J_0 -(z+z^{-1}))^{-1}) \to \det_2 (1+\delta J (J_0 -
(z+z^{-1}))^{-1})
\end{equation}
Clearly, $T(z;\ti J_n)\to T(z;J)$. Thus, $g_n$ converges to the Jost function for
$J$.
\end{proof}

The point of this theorem is that we establish the validity of the GC equations for
defining $u$ in the general context of Theorem~\ref{T5.1}. There is a second point ---
we want to turn this analysis around and directly use the GC equations to prove that,
when $(\alpha)$--$(\gamma)$ hold, $c_n (z)$ and $g_n(z)$ have limits for $z\in\bbD$,
thereby providing a third proof of the hard part of Theorem~\ref{T5.1}. The key is
the following theorem of Coffman \cite{Coff}:

\begin{theorem}[\cite{Coff}]\lb{T7.2} Let $J$ be a $d\times d$ diagonal matrix with
entries $\lambda_1, \dots, \lambda_d$ along the diagonal. Let $A_n$ be a sequence of
$d\times d$ matrices so
\begin{SL}
\item[{\rm{(i)}}]
\begin{equation} \lb{7.9}
\sum_{n=1}^\infty \, \|A_n\|^2 <\infty
\end{equation}

\item[{\rm{(ii)}}] $J$ and $\{J+A_n\}_{n=1}^\infty$ are all invertible.
\end{SL}

Consider solutions $y_n\in\bbC^d$ of
\begin{equation} \lb{7.10}
y_{n+1} = (J+A_n)y_n
\end{equation}
with some initial condition $y_1$. Suppose $\lambda_j$ is a simple eigenvalue with
$\abs{\lambda_j}\neq \abs{\lambda_\ell}$ for $\ell\neq j$. Let
\begin{equation} \lb{7.10a}
f(n) = \prod_{m=1}^{n-1} \, [\lambda_j + (A_m)_{jj}]
\end{equation}
Then there exists an initial condition $y_1$ so that
\begin{equation} \lb{7.11}
\lim_{n\to\infty}\, \f{y_{n,j}}{f(n)}
\end{equation}
exists and is nonzero, while for $\ell\neq j$,
\begin{equation} \lb{7.12}
\lim_{n\to\infty} \, \f{y_{n,\ell}}{f(n)} =0
\end{equation}
\end{theorem}

\begin{remarks} 1. Coffman's result is a discrete analog of continuum (ODE) results of
Hartman-Wintner \cite{HW}. Related work includes Ford \cite{Ford}, Benzaid-Lutz
\cite{BL}, and Janas-Moszy\'nski \cite{JM}.

\smallskip
2. Coffman \cite{Coff} only requires that $\lambda_j$ be a simple eigenvalue and allows
others can have Jordan blocks. In \eqref{7.9}, he allows $2$ to be replaced by $p\in
[1,2]$, but such assumptions imply \eqref{7.9}!

\smallskip
3. A pedagogical presentation of Theorems~\ref{T7.1} and \ref{T7.2} will appear in the
second edition of \cite{OPUC1}. Until that second edition appears, the section will be
available on the web at http://www.math.caltech.edu/opuc.html.
\end{remarks}

\begin{corollary}\lb{C7.3} Let $J$ be a $d\times d$ matrix with simple eigenvalues
$\lambda_1 =1$, $\lambda_2, \dots, \lambda_d$ with $\abs{\lambda_i} \neq
\abs{\lambda_j}$ if $i\neq j$ and $\lambda_2, \dots, \lambda_d\in\bbD$.
Let $y^{(0)}$ be the eigenvector of $J$ with $Jy^{(0)} = y^{(0)}$. Suppose
$A_n$ obeys \eqref{7.9} and
\[
\lim_{N\to\infty}\, \sum_{n=1}^N A_n
\]
exists. Then for any initial condition $y_1$, the solution of \eqref{7.10} obeys
\begin{equation} \lb{7.13}
\lim_{n\to\infty} \, y_n = c(y_1) y^{(0)}
\end{equation}
\end{corollary}

As a preliminary, note that if $\sum_{n=1}^N a_n$ has a limit and
\begin{equation} \lb{7.14x}
\sum_{n=1}^\infty \, \abs{a_n}^2 <\infty
\end{equation}
then $\prod_{j=1}^N (1+a_n)$ has a finite limit which is nonzero if all $a_n \neq -1$.
For $\sum_{n=1}^\infty \log (1+a_n)-a_n$ is absolutely convergent by \eqref{7.14x}.

\begin{proof} By this remark and Theorem~\ref{T7.2}, there are solutions $y_n^{(k)}$
with
\[
y_n^{(k)} \lambda_k^{-n} \to \text{ multiple of eigenvector of $J$ with eigenvalue
$\lambda_k$}
\]
\eqref{7.13} follows since $\lambda_1 = 1$ while $\abs{\lambda_k} <1$ for $k\neq 1$.
\end{proof}

\begin{remark} By using Perron's theorem, one can show that only $\abs{\lambda_j} <1$
for $j\geq 2$ is needed, not $\abs{\lambda_j} \neq \abs{\lambda_k}$.
\end{remark}

Here is the promised third proof of the hard part of Theorem~\ref{T5.1}:

\begin{theorem} \lb{T7.4} Let conditions $(\alpha)$--$(\gamma)$ of Theorem~\ref{T5.1}
hold. Let $c_n, g_n$ be defined by \eqref{7.1}--\eqref{7.3}. Then
\begin{equation} \lb{7.14}
\binom{c_n}{g_n} (z) \to f(z) \binom{(1-z^2)^{-1}}{1}
\end{equation}
\end{theorem}

\begin{proof} Let $J=\left(\begin{smallmatrix} z^2 & 1 \\ 0 & 1 \end{smallmatrix}\right)$
with $z\in\bbD$. $J$ has eigenvalues $z^2 \in\bbD$ and $1$, and the eigenvector for
eigenvalue $1$ is $((1-z^2)^{-1} 1)^t$. Let
\[
A_n = \begin{pmatrix}
-b_{n+1}z & 0 \\
(1-a_{n+1}^2) z^2 - b_{n+1}z & 0
\end{pmatrix}
\]
which obeys the hypothesis of Corollary~\ref{C7.3} by $(\alpha)$--$(\gamma)$. This
corollary plus existence of the limit $\prod_{j=1}^\infty a_j$ imply \eqref{7.14}.
\end{proof}

\begin{remark} One can also apply Theorem~\ref{T7.2} directly to the recursion relation
\eqref{2.1} to see that there exist solutions $\sim Cz^n$, that is, Jost solutions.
\end{remark}

\section{$L^2$ Convergence on the Boundary} \lb{s8}

Our goal in this section is to prove:

\begin{theorem} \lb{T8.1} Let $d\rho$ have the form \eqref{1.1} and suppose
$\{a_n,b_n\}_{n=1}^\infty$ obey $(\alpha)$--$(\gamma)$ of Theorem~\ref{T5.1}. Then
\begin{equation} \lb{8.1}
\lim_{n\to\infty}\, \int_{-2}^2 \abs{p_n(x) -(\sin\theta)^{-1} \Ima (\bar u(e^{i\theta})
e^{i(n+1)\theta})}^2 \, f(x)\, dx =0
\end{equation}
with $\theta =\arccos (\f{x}{2})$, and
\begin{equation} \lb{8.2}
\lim_{n\to\infty} \int \abs{p_n(x)}^2 \, d\rho_\s (x) =0
\end{equation}
\end{theorem}

\begin{remark} Unfortunately, there are some errors in the analogous formula in \cite{OPUC2},
namely, (13.3.15) should have
\begin{equation} \lb{8.2a}
\f{\bar u(x) e^{i(n+1)\theta} - u(x) e^{-i (n+1)\theta}}{2i\sin\theta}
\end{equation}
where it has
\begin{equation} \lb{8.2b}
\f{\bar u(x) e^{i(n-1)\theta} -u(x) e^{i(n-1)\theta}}{4\sin\theta}
\end{equation}
and $p_n$, not $P_n$. As a check, when $b_n\equiv 0$, $a_n\equiv 1$, $p_n (2\cos\theta) =
\f{\sin (n+1)\theta}{\sin\theta}$ and $u\equiv 1$.
\end{remark}

This is an analog of what Szeg\H{o} proved in \cite{Sz20} for OPUC and then translated
\cite{Sz22} to exactly this form for OPRL with $\supp (d\rho)\subset [-2,2]$.
Peherstorfer-Yudistkii \cite{PY} proved precisely this when \eqref{1.5} and \eqref{1.13}
hold. While the underlying core idea behind the proof we use is that of those authors,
our technicalities are much more complex.

For all these proofs, the key is to prove what is essentially a weak $L^2$ convergence
that in the current context is
\begin{equation} \lb{8.3}
\lim_{n\to\infty} \int_0^{2\pi} e^{in\theta} p_n (2\cos\theta) u(e^{i\theta})^{-1}
(1-e^{2i\theta}) \, \f{d\theta}{2\pi} =1
\end{equation}
In Szeg\H{o} case, $u(z)^{-1}$ is an $H^2$-function, so the left side of \eqref{8.3} is
just $\left. [z^n p_n (z+\f{1}{z}) u(z)^{-1} (1-z^2)]\right|_{z=0}$, which converges to
$1$ by the asymptotic result inside the circle. If there are finitely many bound states,
$u(z)^{-1}$ has finitely many poles. Using the fact that eigenfunctions go to zero, it
is easy to accommodate the poles. For the case that Peherstorfer-Yuditskii study, the
argument is more subtle but, by cutting off Blaschke products, still involves a contour
integral around the whole unit circle.

In contrast, our $u(z)^{-1}$ is so singular at $\pm 1$ that we do not see how to directly
deal with the integral in \eqref{8.1}. Instead, we will deal with arcs by mapping a sector
to the unit disk and relating this to distributional convergence of suitable boundary values
of analytic functions. As noted in the introduction, this argument has some elements in
common with work of Denisov-Kupin \cite{DeK} which was done subsequently to our work.
In turn, we were all motivated by some arguments of Killip \cite{Kil}.

The technical core of our proof is the following:

\begin{proposition} \lb{P8.2} Let $(\alpha)$--$(\gamma)$ of Theorem~\ref{T5.1} hold. Fix
a sector on $\bbD$,
\begin{equation} \lb{8.3a}
S=\{z\mid\abs{z}<1, \, 0<\theta_0 \leq \arg z \leq \theta_1<\pi\}
\end{equation}
Then there exist $N$ and $C$ so that for all $n$ and all $z\in S$,
\begin{equation} \lb{8.4}
\biggl| z^n p_n \biggl( z+\f{1}{z}\biggr) (1-z^2) u(z)^{-1} \biggr| \leq
C(1-\abs{z})^{-N}
\end{equation}
Moreover, $C$ and $N$ are uniformly bounded for $S$ fixed for all $\{a_n, b_n\}_{n=1}^\infty$
with
\begin{equation} \lb{8.5}
\sup_N \biggl( \biggl| \, \sum_{n=1}^N \log (a_n)\biggr| + \biggl|\, \sum_{n=1}^N b_n\biggr| \biggr)
+ \sum_{n=1}^\infty \, \abs{a_n -1}^2 + \abs{b_n}^2 <K
\end{equation}
\end{proposition}

\begin{remarks} 1. What is critical is the $C$ and $N$ are $n$ independent. $N$ is also $S$
independent, but $C$ is $S$ dependent and diverges as $S$ approaches the real axis, that
is, as we approach the singularities at $z=\pm 1$.

\smallskip
2. As noted in the introduction, by using ideas of Denisov-Kupin \cite{DeK}, we can likely
prove this with $N=1$; we will have $N=\f52$.

\smallskip
3. We defer the proof until the end of the section.
\end{remarks}

\begin{definition} If $f(z)$ is a function analytic on $\bbD$, with
\begin{equation} \lb{8.6}
f(z) = \sum_{n=0}^\infty \hat f(n) z^n
\end{equation}
and if $f$ obeys
\begin{equation} \lb{8.7}
\abs{f(z)} \leq C (1-\abs{z})^{-N}
\end{equation}
so
\begin{equation} \lb{8.8}
\abs{\hat f(n)} \leq 4 C(n+1)^N
\end{equation}
(obtained by writing $\hat f(n)$ as a contour integral over a circle of radius $1-(n+1)^{-1}$),
we define the distributional boundary values of $f$ by
\begin{equation} \lb{8.9}
\int f(e^{i\theta}) g(e^{i\theta}) \, \f{d\theta}{2\pi} \equiv \sum_{n=0}^\infty
\hat f(n) \int e^{in\theta} g(e^{i\theta})\, \f{d\theta}{2\pi}
\end{equation}
for $C^\infty$ functions $g$ on $\partial\bbD$.
\end{definition}

Power bounds like \eqref{8.4} are important because of

\begin{proposition} \lb{P8.3} Let $f_n(z)$ be a sequence of functions analytic on $\bbD$
so that for some fixed $C,N$ and all $n$,
\begin{equation} \lb{8.10}
\abs{f_n(z)} \leq C(1-\abs{z})^{-N}
\end{equation}
Suppose $f_n\to f_\infty$ uniformly on compacts of $\bbD$. Then the distributional
boundary values converge in {\rm{(}}weak{\rm{)}} distributional sense.
\end{proposition}

\begin{proof} Let $f\in C^\infty (\partial\bbD)$ and $\hat g (-n)$ the integral in
\eqref{8.9}. Write
\[
\sum_{k=0}^\infty \abs{\hat f_n(k) - \hat f_\infty (k)}\, \abs{\hat g(-k)} \leq
\boxed{1} + \boxed{2}
\]
where $\boxed{2} =\sum_{k=K+1}^\infty$ and $\boxed{1}=\sum_{k=1}^K$. Then, by \eqref{8.8}
and \eqref{8.10}, for $g$ fixed, we can choose $K$ so $\abs{\boxed{2}}<\veps$. By the
assumed convergence, $\hat f_n (k) \to \hat f_\infty (k)$ for each $k$, so $\boxed{1} \to
0$.
\end{proof}

\begin{proposition} \lb{P8.4} For any sector $S$ of the form \eqref{8.3a}, there is
an analytic bijection $\varphi:\bbD\to S$ and constant $C$ so that
\begin{equation} \lb{8.11}
(1-\abs{\varphi(z)})^{-1} \leq C (1-\abs{z})^{-1}
\end{equation}
\end{proposition}

\begin{proof} By compactness, we only need to prove this near points where
$\abs{\varphi(z)} =1$. Such points are on the part of $\partial\bbD$ that maps into
$\partial S\cap\partial\bbD$, that is, an arc. At interior points, $\varphi$ is locally
linear and \eqref{8.11} holds, so we only need to worry about neighborhoods, $N$\!, of
the points that map to corners. In suitable local coordinates, $\zeta$, the corner
maps to $0$, $N\cap\bbD$ maps to $\bbC_+\cap \{\zeta\mid\abs{\zeta}<\veps\}$ and
$\varphi$ is transformed to $\varphi(z(\zeta))=\sqrt{\zeta}$ mapping $\bbC_+$ to a
$90^\circ$ corner. In these local coordinates, $1-\abs{\varphi(z)}\sim \Ima \sqrt{\zeta}$
and $1-\abs{z}\sim\Ima\zeta$, and \eqref{8.11} says
\[
\f{1}{\Ima\sqrt{\zeta}} = \f{2\Real\sqrt{\zeta}}{\Ima\zeta} \leq \f{C}{\Ima\zeta}
\]
which is immediate since $\abs{\zeta}$ is small. Thus, \eqref{8.11} holds locally,
and so globally.
\end{proof}

Given a sector $S$ of the form \eqref{8.3a} and $\varphi:\bbD\to S$, we can define
analytic functions on $\bbD$,
\begin{equation} \lb{8.12}
f_n(z) = g_n(\varphi(z))
\end{equation}
where
\begin{equation} \lb{8.13}
g_n(z) = z^n p_n \biggl( z+\f{1}{z}\biggr) u(z)^{-1} (1-z^2)
\end{equation}
This allows us to consider boundary values of $f_n$, and so $g_n$, as distributions.
But $g_n$, and so $f_n$, also has pointwise boundary values, and we want to prove that
the distributional boundary value is given by the function. We have

\begin{lemma} \lb{L8.5} For any sector $S$ of the form \eqref{8.3a}, there is an
$\bbH^2$ function $h$ and a function $q$ analytic in a neighborhood of $\bbD\cup
\bar S$ so that
\begin{equation} \lb{8.14}
u(z)^{-1} = h(z) q(z)
\end{equation}
\end{lemma}

\begin{remark} The point is that $q(z)$ is analytic in $\partial S$, so the boundary
values of $u^{-1}$ on $\partial S$ are given by the well-studied theory of boundary
values of $\bbH^2$ functions \cite{Rudin,Duren}.
\end{remark}

\begin{proof} We use the representation \eqref{5.29}. Define
\begin{equation} \lb{8.15}
H(\theta) = \begin{cases}
\Ima M(e^{i\theta}) & e^{i\theta} \in T\cup\bar T \\
\abs{\sin\theta} & e^{i\theta}\notin T\cup\bar T
\end{cases}
\end{equation}
where $T$ is a slightly enlarged $\partial S\cap\partial\bbD$, but not so
enlarged that it includes $+1$ or $-1$. Let
\begin{equation} \lb{8.16}
h(z) =\exp \biggl( \f{1}{4\pi} \int \f{e^{i\theta} +z}{e^{i\theta}-z}\,
\log \biggl( \f{H(\theta)}{\sin\theta}\biggr)\, d\theta\biggr)
\end{equation}
Since $\f{H(\theta)}{\sin \theta} >0$ (by \eqref{5.30}),
\begin{equation} \lb{8.17}
\int \log \biggl( \f{H(\theta)}{\sin\theta}\biggr) >-\infty
\end{equation}
and, by \eqref{5.30} again,
\begin{equation} \lb{8.18}
\int \biggl( \f{H(\theta)}{\sin\theta}\biggr)\, \f{d\theta}{2\pi} <\infty
\end{equation}
the function \eqref{8.16} is in $\bbH^2$ by the standard approximant argument of
Szeg\H{o} (see \cite[Section~2.4]{OPUC1}).

Thus, we define $q(z)=u(z)^{-1} h(z)^{-1}$ and need to show that it is analytic in a
neighborhood of $\bbD\cup S$. By \eqref{5.29}, we can write
\begin{equation} \lb{8.19}
q(z) =q_1(z) q_2(z) q_3(z) q_4(z)
\end{equation}
where
\begin{equation} \lb{8.20}
q_1(z) = \biggl(\, \prod_{j=1}^\infty a_j\biggr)^{\alpha(a)} \exp \biggl(\tfrac12 \,
\beta(z) \sum_{j=1}^\infty b_j \biggr)
\end{equation}
is clearly analytic away from $\pm 1$. We have that $q_2$ is the inverse of the
renormalized Blaschke product is analytic away from $\bbR$ by Theorem~\ref{T3.3}. Next,
\begin{equation} \lb{8.21}
q_3(z) =\exp\biggl( \f{1}{4\pi} \int Q(z,e^{i\theta}) \log \biggl( \f{\Ima
M}{H(\theta)}\biggr)\, d\theta \biggr)
\end{equation}
is analytic since $\log (\f{\Ima M}{H(\theta)})$ is supported on $\partial\bbD \backslash
(T\cup\bar T)$ and $Q(z,e^{i\theta})$ has singularities only at $z=\pm 1, \pm
e^{i\theta}$. Finally,
\begin{equation} \lb{8.22}
q_4 (z) = \exp \biggl( \f{1}{4\pi} \int [Q(z,e^{i\theta}) - P(z, e^{i\theta})] \log
\biggl( \f{H(\theta)}{\sin\theta}\biggr)\, d\theta\biggr)
\end{equation}
is analytic since $\log (\f{H(\theta)}{\sin\theta})$ is even, so we can replace
$Q-P$ by $Q(z,e^{i\theta}) - \f12 P(z,e^{i\theta}) - \f12 P(z,e^{i\theta})$ and this
kernel is only singular at $z=\pm 1$.
\end{proof}

Given $\theta_0\in [0,2\pi)$, let $R_{\theta_0}$ be the region
\begin{equation} \lb{8.23}
R_{\theta_0} = \biggl\{z\biggm| 1 > \abs{z}  >\tfrac12, \,
\arg (1-e^{-i\theta_0}z) < \f{\pi}{4}\biggr\}
\end{equation}
a region of nontangential approach to $\theta_0$. Define the maximal function,
\begin{equation} \lb{8.24}
M(\theta_0) = \sup_{z\in R_{\theta_0}}\, \abs{u(z)^{-1}}
\end{equation}
Given Lemma~\ref{L8.5}, standard $\bbH^2$ theory \cite{Rudin,Duren,Katz} implies that

\begin{proposition} \lb{P8.6} $u(re^{i\theta})^{-1}$ has boundary values as
$r\uparrow 1$ for a.e.\ $\theta\in (0,2\pi)$. Indeed, for a.e.\ $\theta$,
\begin{equation} \lb{8.25}
\lim_{\substack{ \abs{z}\uparrow 1 \\ z\in R_{\theta_0}}}\, u(z)^{-1} =
u(e^{i\theta_0})^{-1}
\end{equation}
Moreover, for every $\eta >0$,
\begin{equation} \lb{8.26}
\int_\eta^{\pi -\eta} M(\theta_0)^2 \, \f{d\theta}{2\pi} <\infty
\end{equation}
\end{proposition}

This implies

\begin{proposition} \lb{P8.7} Let $S$ be a sector of the form \eqref{8.3a} and $f_n$ be
given by \eqref{8.12}/\eqref{8.13}. Let $\ti S\subset\partial\bbD$ be the image of
$\partial S\cap\partial\bbD$ under $\varphi^{-1}$ and let $T_{f_n}$ denote the
distribution induced by $f_n$ on $\partial\bbD$. Let $t(\theta)$ be a function in
$C^\infty (\ti S^{\intt})$. Then
\begin{equation} \lb{8.27}
T_{f_n}(t) = \int_{\ti S} t(\theta) g_n (\varphi (e^{i\theta}))\, \f{d\theta}{2\pi}
\end{equation}
where $g_n$ is defined by the pointwise boundary value of $u^{-1}$.
\end{proposition}

\begin{proof} By the definition of \eqref{8.9} and the absolute convergence of the
sum,
\begin{equation} \lb{8.28x}
T_{f_n}(t) =\lim_{r\uparrow 1} \int_{\ti S} t(\theta) g_n (\varphi (re^{i\theta})) \,
\f{d\theta}{2\pi}
\end{equation}
By the continuity of $\varphi\restriction \{re^{i\eta}\mid \eta\in\supp(t)\}$, for all
$r$ close enough to $1$, $\varphi (re^{i\theta})\in R_{\varphi(e^{i\theta})}$, so
$g_n (\varphi (re^{i\theta}))\to g_n (\varphi (e^{i\theta}))$ by \eqref{8.25}, and
by \eqref{8.26} and $\abs{t g_n}\leq 2 \abs{t} M \sup_{\abs{t}\leq 1} \abs{z^n p_n
(z+\f{1}{z})}$, we have domination by a function in $L^2$, and so in $L^1$. Thus,
\eqref{8.27} follows from the dominated convergence theorem.
\end{proof}

\begin{proposition} \lb{P8.7A} For each $C^\infty$ function, let $S$ be a sector of
the form \eqref{8.3a} and $\varphi$ an analytic map of $\bbD$ to $S$. Let $t$ be a
$C^\infty$ function with support in $S^\intt$. Then
\begin{equation} \lb{8.27b}
\lim_{n\to\infty}\, \int_{\ti S} t(\theta) g_n (\varphi (e^{i\theta})) \,
\f{d\theta}{2\pi} = \int t(\theta) \, \f{d\theta}{2\pi}
\end{equation}
\end{proposition}

\begin{proof} By \eqref{8.4} and \eqref{8.11}, plus Theorem~\ref{T5.1} (which implies
$g_n(z)\to 1$ for $z\in\bbD$ and so for $z\in S$) and Proposition~\ref{P8.3},
$T_{f_n}\to 1$ in distributional sense. This implies \eqref{8.27b}, given
Proposition~\ref{P8.7}.
\end{proof}

These lengthly preliminaries imply the key to $L^2$ convergence on the boundary:

\begin{theorem} \lb{T8.8} Let $u(e^{i\theta})^{-1}$ be the boundary values of
$u^{-1}$ on $\partial\bbD$. Let $g_n$ be given by \eqref{8.13} on $\partial\bbD$.
Then
\begin{SL}
\item[{\rm{(1)}}]
\begin{equation} \lb{8.27a}
\int_0^{2\pi} \abs{g_n (e^{i\theta})}^2\, \f{d\theta}{2\pi} \leq 2
\end{equation}
\item[{\rm{(2)}}] $g_n\to 1$ in weak-$L^2 (\partial\bbD, \f{d\theta}{2\pi})$.
\end{SL}
\end{theorem}

\begin{proof} (1) We begin with some preliminaries concerning the measure $d\mu_\ac$
on $\partial\bbD$ obtained by using $\theta =\arccos (\f{x}{2})$ to move the a.c.\
part of $d\rho$, that is, $f(x)\,dx$ to $\partial\bbD$. Since $e^{i\theta}\to 2\cos\theta$
is $2-1$ from $\partial\bbD$ to $[-2,2]$ (see \cite[Section~13.15]{OPUC2}),
\begin{equation} \lb{8.28}
d\mu_\ac (\theta) = \abs{\sin\theta} f(2\cos\theta)\, d\theta
\end{equation}
By \eqref{5.2} and standard theory of Stieltjes transforms,
\begin{equation} \lb{8.29}
f(2\cos\theta) = \f{\abs{\Ima M (e^{i\theta})}}{\pi}
\end{equation}
so, by \eqref{5.30},
\begin{align}
d\mu_\ac &= \f{\sin^2 (\theta)}{\abs{u(e^{i\theta})}^2} \, \f{d\theta}{\pi} \lb{8.30} \\
& = \f12\, \f{\abs{1-e^{2i\theta}}^2}{\abs{u(e^{i\theta})}^2} \, \f{d\theta}{2\pi} \lb{8.31}
\end{align}
Thus,
\begin{align*}
\int_0^{2\pi} \abs{g_n (e^{i\theta})}^2 \, \f{d\theta}{2\pi}
&= 2 \int_0^{2\pi}  \abs{p_n (2\cos\theta)}^2 \, d\mu_\ac (\theta) \\
&= 2\int_{-2}^2 \, \abs{p_n(x)}^2 \, d\rho_\ac (x) \\
&\leq 2 \int_{-2}^2 \, \abs{p_n (x)}^2 \, d\rho(x) =2
\end{align*}
proving \eqref{8.27a}.

\smallskip
(2) By (1), the functions $g_n$ are uniformly bounded in $L^2$, so it suffices to prove
that
\begin{equation} \lb{8.32x}
\int t(e^{i\theta}) g_n (e^{i\theta})\, \f{d\theta}{2\pi} \to \int
t(e^{i\theta}) \, \f{d\theta}{2\pi}
\end{equation}
for a total set of $t$'s. If $t$ is $C^\infty$ and supported in some sector $S$ of the
form \eqref{8.3a}, \eqref{8.32x} follows from Proposition~\ref{P8.7A} (there is a
Jacobian to go from $d\varphi(\theta)$ to $d\theta$, but it is $C^\infty$ on $S^\intt$
and occurs on both sides of \eqref{8.32x}). Since such $t$'s are total, (2) is proven.
\end{proof}

\begin{proof}[Proof of Theorem~\ref{T8.1}] Define in $L^2 ([-2,2], f(x)\, dx)$,
\begin{align}
j_n^+(x) &= (2\sin\theta)^{-1} \, \ol{u(e^{i\theta})}\, e^{i(n+1)\theta} \lb{8.32} \\
j_n^-(x) &= \ol{j_n^+(x)} \lb{8.33}
\end{align}
where $\theta(x)\in (0,\pi)$ is given by $x=2\cos (\theta(x))$.
By \eqref{8.29} and \eqref{5.30},
\begin{equation} \lb{8.34}
f(x) = \f{\sin\theta}{\pi \abs{u(e^{i\theta})}^2}
\end{equation}

Thus, by a change of variables,
\begin{align}
\int_{-2}^2 \, \abs{j_n^+(x)}^2 f(x)\, dx  &= \int_{-2}^2 \f14\, \sin^{-2} \theta
\abs{u}^2 \, \f{\sin\theta}{\pi \abs{u}^2}\, dx \notag \\
&= \f12 \int_0^\pi \f{d\theta}{\pi} =\f12 \lb{8.35}
\end{align}

On the other hand, by the same change of variables,
\begin{align}
\langle j_n^-, p_n \rangle_{L^2 (f\,dx)} &= \int_{-2}^2 (2\sin\theta)^{-1}\,
\ol{u(e^{i\theta})}\, e^{in\theta} e^{i\theta}
p_n (2\cos\theta) \, \f{\sin\theta}{\pi \abs{u}^2} \, dx \notag \\
&= \int_0^\pi u(e^{i\theta})^{-1}
\left. \biggl[ (1-z^2) z^n p_n \biggl( z+\f{1}{z}\biggr)\right|_{z=e^{i\theta}}\biggr] \,
\f{-1}{2i}\, \f{d\theta}{\pi} \lb{8.36}  \\
&\to \f{i}{2} \lb{8.37}
\end{align}
by Theorem~\ref{T8.8}. \eqref{8.36} uses
\[
e^{i\theta}\, \f{1}{2\sin\theta} \, \sin\theta\, d(2\cos\theta) = e^{i\theta} \sin\theta
= \f{-1}{2i}\, (1-e^{2i\theta})
\]

Similarly, since $p_n$ is real,
\begin{equation} \lb{8.36a}
\langle j_n^+, p_n \rangle_{L^2 (f\,dx)} \to -\f{i}{2}
\end{equation}
Finally, by the same change of variables that led to \eqref{8.35},
\begin{equation} \lb{8.36b}
\langle j_n^+, j_n^-\rangle_{L^2 (f\,dx)} = \f12 \int_0^\pi e^{-2i (n+1)\theta} \,
\f{u(e^{i\theta})}{\,\ol{u(e^{i\theta})}\,}\, \f{d\theta}{\pi} \to 0
 \end{equation}
since $\f{u}{\bar u}\in L^2 (\partial\bbD, \f{d\theta}{2\pi})$ and
$e^{-2i(n+1)\theta)}\to 0$ weakly.

Now,
\begin{equation} \lb{8.38}
p_n(x) -(\sin\theta)^{-1} \Ima (\bar u e^{i(n+1)\theta}) =
p_n(x) - i^{-1} [j_n^+ - j_n^-]
\end{equation}
so, by \eqref{8.35}, \eqref{8.36}, \eqref{8.36a}, and \eqref{8.37} (all norms in $L^2
(f\,dx)$),
\begin{align*}
0 &\leq \liminf \|\text{LHS of \eqref{8.38}}\|^2 \\
&=\liminf (\|p_n\|^2 + \|j_n^+\|^2 + \|j_n^-\|^2 - 2\Real \langle j_n^+, j_n^-\rangle \\
& \qquad \qquad \qquad + 2\Real \langle p_n, ij_n^+\rangle - 2\Real \langle p_n, ij_n^-\rangle ) \\
&= \liminf (\|p_n\|^2 + \tfrac12 + \tfrac12 - 0 - 1-1) \\
&= \liminf \|p_n\|^2 -1
\end{align*}

Thus, since $\|p_n\|^2_{L^2(f\, dx)} \leq 1$, we conclude
\[
\lim\|p_n\|_{L^2(f\,dx)}^2 =1
\]
so $\|p_n\|_{L^2 (d\rho_\s)}\to 0$ and, by the above calculation, $\text{LHS of \eqref{8.38}} \to 0$.
\end{proof}

Thus, Theorem~\ref{T8.1} is reduced to the proof of Proposition~\ref{P8.2}, to which we now turn.
As a preliminary, we want to exploit the proof of Lemma~\ref{L8.5}:

\begin{proposition} \lb{P8.9} Let $S$ be a sector of the form \eqref{8.3a} and let $K$ be given.
Then there exists $C$ depending only on $S$ and $K$ so that if \eqref{8.5} holds, then for
$z\in S$,
\begin{equation} \lb{8.40}
\abs{u(z)^{-1}} \leq C(1-\abs{z})^{-1/2}
\end{equation}
\end{proposition}

\begin{proof} As in Lemma~\ref{L8.5}, we construct the factorization \eqref{8.14}. The proof
shows that $\|h\|_{\bbH^2}$ and $\|q\|_{S,\infty} \equiv \sup_{z\in S} \abs{q(z)}$ are
bounded by constants $C_1$ and $C_2$ depending only on $K$ and $S$.

Let $h(z)=\sum_{j=0}^\infty h_j z^j$, then $\|h\|_{\bbH^2} = (\sum_{j=0}^\infty \abs{h_j}^2)^{1/2}$,
so, by the Schwarz inequality,
\begin{align*}
\abs{h(z)} &\leq \|h\|_{\bbH^2} \biggl( \, \sum_{j=0}^\infty \, \abs{z}^{2j}\biggr)^{1/2} \\
&= \|h\|_{\bbH^2} (1-\abs{z}^2)^{-1/2} \\
&\leq C_1 (1-\abs{z})^{-1/2}
\end{align*}
so, by \eqref{8.14}, we have \eqref{8.40} with $C=C_1 C_2$.
\end{proof}

Our proof will exploit \eqref{2.10} where $y_n = \ti w_n$ and $x_n = c_{n-1}$. We are
interested in controlling $u^{-1} x_n$, which means controlling $u^{-1} y_n^{-1}$ and
ratios $y_n/y_{n-j}$, that is, the functions $u^{-1} \ti w_n^{-1}$ and $\ti w_n/\ti w_{n-j}$.
So we turn first to $u^{-1} \ti w_n^{-1}$ and then $\ti w_n/\ti w_{n-j}$.

Let $J^{(n)}$ be the Jacobi matrix given after \eqref{5.2} and make the $J$-dependence of
$w_n$ explicit. Then:

\begin{proposition} \lb{P8.10} Let $\ti w_n (z,J)$ be given by \eqref{2.2}/\eqref{2.3}. Then
\begin{equation} \lb{8.41}
\ti w_n (z,J)^{-1} u(z,J)^{-1} = a_n u(z,J^{(n)})^{-1}
\end{equation}
In particular, if \eqref{8.5} holds, then for any $S$ obeying \eqref{8.3a}, there is a
$C$ so
\begin{equation} \lb{8.42}
\sup_n \, \abs{\ti w_n (z,J)^{-1} u(z,J)^{-1}} \leq C (1-\abs{z})^{-1/2}
\end{equation}
\end{proposition}

\begin{remarks} 1. In order to get \eqref{8.40}, one does not need a bound on
$\sup_N \abs{\sum_1^N b_n}$ but only on $\lim_N \abs{\sum_1^N b_n}$ (and similarly for
$\log (a_n)$). But to get \eqref{8.42}, we need control on $\sup_N \lim_M (\sum_N^{N+M}
b_n)$ --- and that is why we state \eqref{8.5} in the form we do.

\smallskip
2. One can also prove this result using the fact that $uw_n$ is the unique solution
asymptotic to $z^n$.
\end{remarks}

\begin{proof} By \eqref{5.6},
\[
\ti w_n (z,J) = \biggl( \f{M(z)}{z}\biggr) \biggl( \f{a_1 M_1 (z)}{z}\biggr) \dots
\biggl( \f{a_{n-1} M_{n-1}(z)}{z}\biggr)
\]
from which it follows that
\begin{equation} \lb{8.42x}
\ti w_{n+k}(z,J) = a_n \ti w_n (z,J) \ti w_k (z,J^{(n)})
\end{equation}
Taking $k$ to infinity using \eqref{2.3} and $\ti w_\infty = u^{-1}$, we obtain
\[
u(z,J)^{-1} = a_n \ti w_n (z,J) u(z,J^{(n)})^{-1}
\]
which is \eqref{8.41}.
\end{proof}

\begin{proposition} \lb{P8.11} For any $J$\!,
\begin{SL}
\item[{\rm{(i)}}]
\begin{equation} \lb{8.43}
\abs{\ti w_n (z,J)} \leq \f{\pi}{4}\, \abs{z}^{-n} \veps^{-1} (1-\abs{z})^{-1}
\end{equation}
when $\arg z\in (\veps, \pi -\veps)$ and $0 <\abs{z} <1$.
\item[{\rm{(ii)}}]
\begin{equation} \lb{8.43a}
\biggl| \f{\ti w_{n+k}(z,J)}{\ti w_n (z,J)}\biggr| \leq \f{\pi}{4} \, a_n \abs{z}^{-k}
\veps^{-1} (1-\abs{z}^{-1})
\end{equation}
when $\arg z\in (\veps, \pi-\veps)$ and $0<\abs{z} <1$.
\end{SL}
\end{proposition}

\begin{proof} (i) By \eqref{2.2},
\begin{align*}
\abs{\ti w_n (z,J)} &\leq \|(z+z^{-1} -J)^{-1} \| \, \abs{z}^{-n} \\
&\leq \abs{\Ima (z+z^{-1})}^{-1} \, \abs{z}^{-n}
\end{align*}
since $\sigma(J)\subset\bbR$ and $J$ is selfadjoint. But if $z=re^{i\varphi}$ with $0 < r
< 1$,
\begin{equation} \lb{8.44}
\abs{\Ima (z+z^{-1})} = (r^{-1}- r) \abs{\sin\varphi}
\end{equation}
For $\varphi\in (\veps, \pi-\veps)$,
\begin{equation} \lb{8.45}
\abs{\sin\varphi} \geq \f{2}{\pi}\, \veps
\end{equation}
and
\begin{equation} \lb{8.46}
r^{-1} -r = (1-r) (r^{-1} +1) \geq 2(1-r)
\end{equation}
Thus \eqref{8.43} holds.

\smallskip
(ii) \eqref{8.43a} follows immediately from \eqref{8.43} and \eqref{8.42x}.
\end{proof}

\begin{proof}[Proof of Proposition~\ref{P8.2}] By \eqref{2.10} and the
 proof of
Theorem~\ref{T2.2},
\begin{equation} \lb{8.47}
\abs{c_n(z)} \leq \sum_{j=0}^{n-1} a_{n-j}^{-1} \, \f{\ti w_{n+1}}{\ti w_{n+1-j} \ti
w_{n-j}} \, z^{2j} + z^{2n} \, \f{\ti w_{n+1}}{\ti w_1} 
\end{equation}
Define
\[
A=\sup_n\, (\abs{a_n}, \abs{a_n}^{-1}) <\infty
\]
since $a_n\to 1$. Then, by \eqref{8.47},
\begin{equation} \lb{8.48}
\begin{split}
\sup_n\, \abs{c_n(z) u(z)^{-1}} & \leq A(1-\abs{z})^{-1} \sup_{n,j} \,
\biggl| \f{\ti w_{n+1} z^j}{\ti w_{n+1-j}}\biggr| \sup_{n,j} \,
\f{1}{\abs{\ti w_{n-j} u}} \\
&\qquad \quad + 
\sup_n \, \biggl| \f{\ti w_{n+1} z^n}{\ti w_1}\biggr| \abs{u}^{-1}
\end{split}
\end{equation}

By \eqref{8.43a}, with $\veps$ chosen so $S\subset \{z\mid\arg z\in
(\veps, \pi-\veps)\}$, we have
\[
\sup_{n,j}\, \biggl| \f{\ti w_{n+1} z^j}{w_{n+i-j}}\biggr| \leq
\f{\pi A}{4}\, \veps^{-1} (1-\abs{z})^{-1}
\]
By \eqref{8.42},
\[
\sup_{n,j}\, \biggl| \f{1}{\ti w_{n-j} u}\biggr| \leq C(1-\abs{z})^{-1/2}
\]

Thus, \eqref{8.48} implies that
\[
\biggl| z^n p_n \biggl( z+\f{1}{z}\biggr) (1-z^2) u(z)^{-1} \biggr| \leq
C (1-\abs{z})^{-5/2}
\]
where $C$ depends on $S$ and the constant $K$ in \eqref{8.5}.
\end{proof}

This completes the proof of Theorem~\ref{T8.1}. We note that what is missing
is proving that wave operators exist. One can use \eqref{8.1} to prove that for
a.e.\ $\theta$, the transfer matrix is bounded along a subsequence and so obtain
an alternate proof of the Deift-Killip theorem \cite{DK} when $(\beta)$, $(\gamma)$
hold. Of course, we use the Killip-Simon sum rule which also implies that theorem.
Still it is interesting that one has plane waves like solutions in $L^2$ sense in
$[-2,2]$. The deeper and interesting open question is pointwise convergence for
a.e.\ $\theta$.

\section{Bound States} \lb{s9}

One knows that with regard to Szeg\H{o} asymptotics, sometimes simple-looking
assumptions are really quite restrictive: for instance (see, e.g.,
\cite[Chapter~13]{OPUC2}), if $\supp (d\mu)\subset [-2,2]$, then $(\beta)$
of Theorem~\ref{T5.1} implies all of $(\alpha)$--$(\gamma)$ and all the
other hypotheses of that theorem. Also (see \cite{SZ}), if $(\beta)$ holds
and $f$ is given by \eqref{1.1}, then
\begin{equation} \lb{9.1}
\sum_{j=1,\pm}^{N_\pm} \, (\abs{E_j^\pm}-2)^{1/2} < \infty
\end{equation}
if and only if
\begin{equation} \lb{9.2}
\int_{-2}^2 (\log f)(4-x^2)^{-1/2}\, dx > -\infty
\end{equation}
Here we want to show that $(\beta)$, $(\gamma)$ alone do not imply spectral restrictions.
In particular, we want to show that for each $q< \f32$, there is a Jacobi matrix obeying
$(\alpha)$--$(\gamma)$ where
\begin{equation} \lb{9.3}
\sum_{j=1,\pm}^{N_\pm} (\abs{E_j^\pm}-2)^q =\infty
\end{equation}
Of course, by \cite{KS}, $(\alpha)$ implies
\begin{equation} \lb{9.4}
\sum_{j=1,\pm}^{N_\pm} (\abs{E_j^\pm}-2)^{3/2} <\infty
\end{equation}

Our construction will have $a_n\equiv 1$ and $b_n$ nonzero in blocks. In
\cite[Section~13.9]{OPUC2}, examples with $b_n$ nonzero in a sequence of isolated
points are constructed where $(\alpha)$--$(\gamma)$ hold and \eqref{9.3} holds for
$p$ arbitrarily close to $1$. So this section improves that result. Our
construction is closely related to that in Theorem~5.12 of \cite{DHKS}.

Pick $\alpha$ in $(\f12,1)$ and $p$ so that
\begin{equation} \lb{9.5}
\f{\alpha}{1-\alpha} > p > \f{\alpha}{2-\alpha}
\end{equation}
We will eventually take $\alpha$ to $\f12$ and $p-\f{\alpha}{2-\alpha}\to 0$. Pick $M_0$
and $C_1$ so for $m\geq M_0$, the distances between the blocks $B_m \equiv [m^{p+1} -C_1
m^p, m^{p+1} + C_1 m^p]$ for $m=M_0, M_0 +1, \dots$ are each at least $2$. This is easy
to do if one fixes $C_1 < \f12 (p+1)$. We should use $[C_1 m^p]$, but for notational
simplicity, we will pretend that $C_1m^p$ is an integer. We pick $b_n$ by
\begin{equation} \lb{9.6}
b_n = \begin{cases}
n^{-\alpha} & n\in B_{2k}, \, 2k\geq M_0 \\
-n^{-\alpha} & n\in B_{2k+1}, \, 2k +1 \geq M_0 \\
0 & \text{otherwise} \end{cases}
\end{equation}

\begin{lemma}\lb{L9.1} $a_n\equiv 1$ and $b_n$ in \eqref{9.6} obey $(\alpha)$--$(\gamma)$
of Theorem~\ref{T5.1}.
\end{lemma}

\begin{proof} Since $\abs{b_n}\leq n^{-\alpha}$ and $\alpha > \f12$, condition $(\alpha)$ holds
and $(\beta)$ is trivial. So we only need to check $(\gamma)$. Since $\f{d}{dn} n^{-\alpha}
 = -\alpha n^{-\alpha-1}$, if $n\in B_m$, then
\begin{align*}
\abs{n^{-\alpha} - m^{-(p+1)\alpha}} &\leq Cm^{-(p+1)(\alpha+1)} m^p \\
&= Cm^{-[(p+1)\alpha+1]}
\end{align*}

Thus,
\begin{equation} \lb{9.7}
\biggl|\, \sum_{n\in B_m} \, \abs{b_n} - 2C_1 m^{-\alpha(p+1) + p}\biggr|
\leq C_2 m^{-1} m^{-\alpha (p+1) + p}
\end{equation}

We claim
\begin{equation} \lb{9.8}
\alpha (p+1) > p
\end{equation}
which implies, first, that the estimate on the right of \eqref{9.7} is absolutely summable
and, second, that $\sum_m (-1)^m m^{-\alpha (p+1) +p}$ is conditionally summable, proving
$(\gamma)$.

To prove \eqref{9.8}, note that it is equivalent to $\alpha > p(1-\alpha)$ or $p<
\f{\alpha}{1-\alpha}$, which is true by \eqref{9.5}.
\end{proof}

For $m$ even, we will pick $\varphi_m$ to be the trial vector supported in $B_m$, which
is $1$ at the center of $B_m$ (i.e., at $n=m^{p+1}$), $0$ at the end points, and constant
slope in between. For $m$ odd, we do the same construction and then multiply by $(-1)^n$.

Consider $m$ even first. Since $a_n\equiv 1$,
\begin{align}
\langle\varphi_m, (J_0 -2)\varphi_m\rangle &= -\sum_n \,
\abs{\varphi_m (n+1) -\varphi_m (n)}^2 \lb{9.9} \\
&\geq - m^p \biggl[ \f{C}{m^p}\biggr]^2 \notag \\
&= -C^2m^{-p} \lb{9.10}
\end{align}
since the slope $\sim m^{-p}$ and there are $O(m^p)$ nonzero terms in the sum
\eqref{9.9}. On the other hand, since $b_n > C_3 m^{-\alpha (p+1)}$ on $B_m$ and, on
average, $\abs{\varphi_m}^2 \geq \f14$ on $B_m$,
\begin{equation} \lb{9.11}
\langle \varphi_m, b\varphi_m\rangle \geq C_3 m^p m^{-\alpha (p+1)}
\end{equation}

It is easy to see that $p>\f{\alpha}{2-\alpha}$ is equivalent to $2p > \alpha (p+1)$, so
for $m\geq M_1$ for some $M_1$,
\begin{equation} \lb{9.12}
\langle \varphi_m (J_0 +b -2)\varphi_m\rangle \geq \tfrac12\, C_3 m^p m^{-\alpha (p+1)}
\end{equation}

Since $\|\varphi_m\|^2 \leq C_4 m^p$, we see that
\begin{equation} \lb{9.13}
\f{\langle \varphi_m, (J_0 + b-2)\varphi_m\rangle}{\|\varphi_m\|^2} \geq C_4 m^{-\alpha
(p+1)}
\end{equation}
for $m$ even. Similarly, for $m$ odd,
\begin{equation} \lb{9.14}
\f{\langle\varphi_m, (J_0 + 2+b)\varphi_m\rangle}{\|\varphi_m\|^2} \leq -C_4
m^{-\alpha(p+1)}
\end{equation}
Since $\langle\varphi_m, \varphi_k\rangle =0=\langle\varphi_m, (J_0 +b)\varphi_k\rangle$
for $m\neq k$, a variational argument proves for $m$ large,
\begin{equation} \lb{9.15}
\abs{\abs{E_m^\pm} -2} \geq \tfrac12\, C_4 m^{-\alpha(p+1)}
\end{equation}

Thus, \eqref{9.3} holds if $q\alpha (p+1) <1$. Taking $\alpha\downarrow \f12$, $p\downarrow
\f13$, we see $q\uparrow ((\f12)(\f43))^{-1} = \f32$. Thus,

\begin{theorem} \lb{T9.1} For any $q<\f32$, there is a set of Jacobi parameters for which
$(\alpha)$--$(\gamma)$ of Theorem~\ref{T5.1} hold, but for which \eqref{9.3} also holds.
\end{theorem}

\section{A Remark on Schr\"odinger Operators} \lb{s10}

In this section, we want to show how the ideas of Section~\ref{s6} provide a
simple proof of

\begin{theorem}\lb{T10.1} Suppose $V\in L^2 (0,\infty)$ and
\begin{equation} \lb{10.1}
\lim_{x\to \infty}\, \int_0^x V(y)\, dy
\end{equation}
exists. Then for any $\kappa$ with $\kappa >0$, there is a solution of
\begin{equation} \lb{10.2}
-u'' + Vu = -\kappa^2 u
\end{equation}
so that
\begin{equation} \lb{10.3}
\lim_{x\to\infty}\, e^{\kappa x} u(x) =1
\end{equation}
\end{theorem}

This result is not new. It was proven by Hartman-Winter \cite{HW} using sophisticated ODE
asymptotic methods. Even with the simplification of Harris-Lutz \cite{HL}, the proof is
involved (see Eastham \cite{East} for a particularly clear discussion of this proof).
Here, as in Section~\ref{s6}, we will use renormalized determinants to construct $u$. The
same argument shows that if \eqref{10.1} does not have a finite limit, then there is a
solution so $u(x)/\exp[f(x)] \to 1$, where
\begin{equation} \lb{10.4}
f(x) = -\kappa x + \f{1}{2\kappa} \int_0^x V(y)\, dy
\end{equation}
also a result of Hartman-Wintner \cite{HW}.

In the argument below, we will use unfactorized kernels (i.e., $VG_0$) rather than
factorized kernels (i.e., $V^{1/2} G_0 \abs{V}^{1/2}$). By using factorized kernels,
one should be able to extend this theorem to the case where $V\in L^2$ is replaced
by $\sum_n (\int_n^{n+1} \abs{V(x)}\, dx)^2 <\infty$.

The starting point is a formula of Jost-Pais \cite{JP} for the Jost function extended
to get the Jost solutions.

\begin{proposition} \lb{P10.2} Let $G_0 (-\kappa^2)$ be the operator $(H_0 +\kappa^2)^{-1}$
where $H_0$ is $-\f{d^2}{dx^2}$ with $u(0)$ boundary conditions, so $G_0$ has integral
kernel
\begin{equation} \lb{10.5}
G_0 (x,y; -\kappa^2) = (2\kappa)^{-1} [e^{-\kappa \abs{x-y}} -e^{-\kappa (x+y)}]
\end{equation}
For any $V\in L^2$ of compact support and any $x_0 >0$, let $K(x_0;\kappa)$ be the
operator with integral kernel
\begin{equation} \lb{10.6}
K(x,y;x_0;\kappa) = V(x+x_0) G_0 (x,y; -\kappa^2)
\end{equation}
Then $K$ is trace class and
\begin{equation} \lb{10.7}
u(x_0) \equiv e^{-\kappa x_0} \det (1+K(x_0;\kappa))
\end{equation}
obeys \eqref{10.2}, and for $x$ large,
\begin{equation} \lb{10.8}
u(x) = e^{-\kappa x}
\end{equation}
\end{proposition}

\begin{remark} \eqref{10.8} is trivial since $V$ has compact support, which means
$K\equiv 0$ for $x_0$ large.
\end{remark}

\begin{proof} This is essentially Proposition~2.9 of \cite{S276}. That paper uses
a factorized kernel, but by the Birman-Solomyak theorem (see \cite[Chapter~4]{TI}, $K$ is
trace class, and so the determinants are equal. As noted, \eqref{10.8} is immediate.
\end{proof}

\begin{proposition} \lb{P10.3} If $V\in L^2$ and \eqref{10.1} holds, then
$K(x_0; \kappa)$ is Hilbert-Schmidt and
\begin{equation} \lb{10.9}
u(x_0) = e^{-\kappa x_0} \det_2 (1+ K( x_0; \kappa)) \exp\biggl( (2\kappa)^{-1}
\int_{x_0}^\infty V(y) [1-e^{-2\kappa y}]\, dy \biggr)
\end{equation}
Moreover, $u$ obeys \eqref{10.2}.
\end{proposition}

\begin{proof} If $V$ has compact support, \eqref{10.9} is just \eqref{10.7} since
$\tr (K(x_0;\kappa))=\int_{x_0}^\infty V(y) (2\kappa)^{-1} [1-e^{-2\kappa y}]\, dy$
and we have \eqref{6.11}. Given general $V$\!, let $V_L(x)$ be given by
\[
V_L(x) = \begin{cases}
V(x) & x\leq L \\
0 & x >L
\end{cases}
\]
and $u_L$ given by \eqref{10.1}. Since $V\in L^2$, $K_L (x_0, \kappa) \to
K(x_0; \kappa)$ in Hilbert-Schmidt norm, so $\det_2$ converges. By \eqref{10.1},
the exponentials converge. Thus, $u_L\to u$. This means $u$ is a distributional
solution of \eqref{10.2} and so, by elliptic regularity, a solution $L^2$ at infinity.
\end{proof}

\begin{proof}[Proof of Theorem~\ref{T10.1}] $K(x_0;\kappa)\to 0$
in Hilbert-Schmidt norm as $x_0\to\infty$, so $\det_2 (1+K(x_0; \kappa))\to 1$. The
integral goes to $0$ as $x_0\to\infty$. Thus, $u(x_0) e^{\kappa x_0} \to 1$.
\end{proof}

The point here is that it is natural to try to construct $u$ as a limit of $u_L$'s,
and then prove asymptotics of $u$. The fact that we have an explicit formula in terms
of renormalized determinants allows us to control both the limit as $L\to\infty$
and then as $x\to\infty$.

\bigskip

\end{document}